\magnification=\magstep1
\baselineskip=18pt
\parskip=6pt
\input amstex
\documentstyle{amsppt}
\vsize=9truein
\loadmsbm
\loadeufm
\loadeufb
\loadeusb
\loadeurm
\loadeurb
\loadeusm
\NoBlackBoxes
\topmatter

\centerline{\bf ON GENERALIZED FOURIER TRANSFORMS}
\medskip
\centerline{\bf FOR STANDARD $L$--FUNCTIONS}
\bigskip
\centerline{\bf with an appendix by Wen--Wei Li}
\vskip.5truein
\centerline{\bf F. SHAHIDI*}
\vskip-.5truein
\footnote"{}"{*Partially supported by the NSF grant DMS--1500759 \hfil\break} 

\loadmsbm
\loadeufm
\UseAMSsymbols
\def\ds{\displaystyle}
\def\ab{\text{ab}}
\def\id{\text{id}}

\def\Hom{\text{Hom}}
\def\der{\text{der}}
\def\det{\text{det}}

\def\Sat{\text{Sat}}
\def\sc{\text{sc}}
\def\on{\overline n}
\def\sym{\text{sym}}
\def\Sym{\text{Sym}}

\def\std{\text{std}}
\def\End{\text{End}}

\def\bC{\Bbb C}
\def\bQ{\Bbb Q}
\def\bA{\Bbb A}
\def\bL{\Bbb L}

\def\bZ{\Bbb Z}

\def\bG{\Bbb G}
\def\bR{\Bbb R}
\def\bF{\Bbb F}
\def\bP{\Bbb P}

\def\Ind{\text{Ind}}

\def\om{\overline m}

\def\oP{\overline P}

\def\oI{\overline I}

\def\oN{\overline N}
\def\sF{\Cal F}

\def\sH{\Cal H}

\def\sJ{\Cal J}

\def\sM{\Cal M}
\def\sR{\Cal R}
\def\sS{\Cal S}
\def\sQ{\Cal Q}

\def\sL{\Cal L}

\def\tw{\widetilde w}
\def\tH{\widetilde H}
\def\tP{\widetilde P}
\def\tN{\widetilde N}
\def\tG{\widetilde G}
\def\tmu{\widetilde \mu}
\def\tpi{\widetilde \pi}
\def\tchi{\widetilde \chi}
\def\tal{\widetilde \alpha}
\def\tbe{\widetilde \beta}
\def\tga{\widetilde \gamma}
\def\tlam{\widetilde \lambda}
\def\hF{\hat F}
\def\hT{\hat T}
\def\hM{\hat M}

\bigskip
\bigskip\noindent
\abstract{Any generalization of the method of Godement--Jacquet on principal $L$--functions for $GL(n)$ to other groups as perceived by Braverman--Kazhdan/Ngo requires a Fourier transform on a space of Schwartz functions. In the case of standard $L$--functions for classical groups, a theory of this nature was developed by Piatetski--Shapiro and Rallis, called the doubling method. It was later that Braverman and Kazhdan, using an algebro--geometric approach, different from doubling method, introduced  a space of Schwartz functions and a Fourier transform, which projected onto those from doubling method. In both methods a normalized intertwining operator played the role of the Fourier transform. The purpose of this paper is to show that the Fourier transform of Braverman--Kazhdan projects onto that of doubling method. In particular, we show that they preserve their corresponding basic functions. The normalizations involved are not the standard ones suggested by Langlands, but rather a singular version of local coefficients of Langlands-Shahidi method. The basic function will require a shift by 1/2 as dictated by doubling construction, reflecting the global theory, and begs explanation when compared with the work of Bouthier--Ngo--Sakellaridis. This matter is further discussed in an appendix by Wen--Wei Li.}
\endabstract
\endtopmatter

\bigskip
\centerline{\bf Introduction}

In a series of papers [BK1,2,3], Braverman and Kazhdan proposed a generalization of the work of Godement and Jacquet [GJ] on principal $L$--functions for $GL(n)$ to an arbitrary reductive group $G$ and a finite dimensional irreducible representation $\rho$ of its $L$--group. This approach has now been taken up by Ngo [BNS,ChN,N1,2,3] and others, who have been defining some of the objects which generalize those in [GJ], for example, a reductive monoid $M_\rho$, studied and classified by E.~B.~Vinberg [V] in characteristic zero, which replaces the simple algebra $M_n$ whose group of units is $GL_n$, where the theory was worked out in [GJ] as pointed out.

Next, one has to define a space of Schwartz functions on the $F$--points of the group of units of the monoid $M_\rho$, $F$ a local field, since in practice this is the space of functions with which one would be working. They are expected to be restrictions of smooth functions of compact support on $M_\rho(F)$. With such a space of functions in hand any generalization of the work in [GJ] demands a Fourier transform acting on this space.

The main purpose of this paper is to address these Fourier transforms in the only case where a theory close to [GJ] exists, namely, the case of standard $L$--functions for classical groups. This is what is usually called ``the doubling method'' introduced by Piatetski--Shapiro and Rallis [GPSR,PSR] in mid-eighties; with contributions from Cogdell in [GPSR]. Its local theory was later detailed, improved and completed by Lapid and Rallis in [LR]. The cases of unitary groups were also  addressed in [L] by Jian--Shu Li who followed the approach in [GPSR], and S.~Yamana [Y1,2] who pursued the approach of Lapid--Rallis [LR]. Gan developed the theory for double covers of symplectic groups in [Gan].

In all these cases, both notions of a Schwartz space, in the form of the space of a parabolically induced representation, and the Fourier transform, as a normalized intertwining operator on this space, are evident, but never explicitly acknowledged. This was made explicit by Braverman and Kazhdan in [BK2,3] through an algebro--geometric approach, and fairly independent of it, as opposed to the doubling method where these objects are all defined by analytic tools, dictated by global theory [GPSR,PSR].

In a recent paper [Li], Wen--Wei Li has carefully analyzed the geometry and other aspects involved in the theory and expanded upon the work in [BK2], but did not address in any detail the Fourier transforms introduced in [BK2] and  their connection with normalized intertwining operators of doubling method; in particular, with the exception of some discussions in Section 8, e.g., Example 8.1.7, the behavior of the basic function under the transform was not discussed in any detail. This has now been supplemented in an appendix by Li [Li2] to the present paper which we will discuss later.

The main result of this paper is to show that the Fourier transforms of Braverman--
--Kazhdan are defined by the same normalizing factors as those in doubling method (Corollary 6.38). Moreover, we show that the Fourier transform preserves the basic function, the unramified function for which the zeta function (5.13) gives the unramified $L$--function precisely, both in doubling construction (Proposition 6.13) and that of Braverman--Kazhdan (Proposition 6.45). We point out that Proposition 6.45, in which basic function is defined by equation (6.16), is more in the spirit of a generalization of Godement--Jacquet, since the zeta function (6.17) is, in fact, an integration of a Schwartz function against a matrix coefficient, but over the $F$--points of the units $M_{ab}\times G$ of the monoid $X$ attached to the standard representation of the classical group $G$. We refer to definitions and discussions in Sections 2, 4 and 6 here, as well as Section 7 of [Li], as well as the appendix [Li2],  concerning these objects. Also see Remark 6.47. 
(Braverman and Kazhdan assumed the group to be simply connected and split, thus covering only the case $G=Sp(n)$ from [GPSR,PSR,LR]. But that is not an issue when discussing normalizing factors.)

By the nature of doubling method, one needs a shift $s-\frac 12$ to obtain the unramified standard $L$--function (Proposition 5.17). In fact, such a shift is imposed upon us by the use of Eisenstein series which enters the method globally [GPSR] and cannot be avoided. In particular, we do not seem to get the shift suggested by Ngo [BNS,N3]. The case of doubling for $GL(n)$ also has this shift [Y2] and is different from Godement--Jacquet [GJ]. Of course, in doubling for $GL(n)$ one gets a product of the principal $L$--function and its dual, rather than just the principal $L$--function as in [GJ]. It should be pointed out that the case of $GL(n)$ is needed to provide us with the local $L$--function for the unitary group over a global field at a split place.

The nature of this shift is addressed and put in context, among other things, in an appendix [Li2] to the present paper by Wen--Wei Li, which came about after the first distribution of this manuscript.

Sections 4, 5, and 6 are devoted to answering these questions. Section 4 is a friendly interpretation of the normalizing factors and Fourier transforms of~Braverman--Kazhdan [BK2]. Section 5 is an exposition on doubling method following [LR]. In particular, we determine explicitly the normalizing factor $\eta(s)$ as a correction factor in terms of the functions $a_H$ and $d_H$ of [PSR] and use them in the definition of our basic function (6.6) in the setting of doubling method. As we explain in Section 5, the work in [LR] introduces $\eta(s)$ at other (ramified) places as a {\bf degenerate local coefficient} $c(s,\chi, A,\psi)^{-1}$, defined the same way as the original local coefficients in [Sh1], but for a degenerate representation of the Levi subgroup, namely, just a character [LR]. We remark that this degenerate local coefficient, and those original ones in [Sh1], normalize intertwining operators by $\gamma$--factors, rather than the way needed in the trace formula and Eisenstein series as suggested by Langlands [A,La,Sh3].

To prove the equality of normalizing factors of Braverman--Kazhdan and those of doubling, we need explicit calculations in terms of adjoint actions as in [Sh1,3] which are done very explicitly in [GPSR] and used in Section 6 here, repeatedly.

The paper starts with a motivational section, Section 1, on what we call a $\gamma$--distribution, a distribution which when evaluated at an irreducible representation, gives the corresponding $\gamma$--factor. We introduce a pair of conjugacy--invariant distributions $\Phi^{GJ}_{\psi}$ and $\Phi^{LS}_{\psi}$, which give the $\gamma$--factor of Godement--Jacquet's principal $L$--functions and Langlands--Shahidi method [Sh1,3], respectively (Proposition 1.17 and paragraph 1.20), pointing similarities of defining methods.

Section 2 gives a review of Vinberg's monoids [V] and the determination of the monoid $M_\rho=M^\lambda$ as a background for Sections 3 and 6, where $\lambda$ is the highest weight of $\rho$, following Ngo [N1,N2]. We conclude Section 2 by comparing this approach to that of conjectural generalization of Langlands--Shahidi method to Kac--Moody groups. We refer to paragraph 2.21 for some common features.

If one restricts the conjectural Fourier transform on the group to a maximal torus, it would behave like a Hankel transform, and this is what is addressed by Ngo in a recent paper [N3]. Our Section 3 is devoted to a short review of his work in [N3], as well as the basic functions and his proposed shift presented there. The equality of our Fourier transforms with the corresponding Hankel transform in [N3] now seems to be close in hand due to  progress made recently in a group activity during the Workshop on ``Beyond Endoscopy and Trace Formula'' at AIM, December 4--8, 2017.

We did not address the global issues where a Poisson summation formula would be needed in this paper at all. But as discussed in [BK2] that may be provided by the functional equation of the Eisenstein series. This is discussed in a recent preprint of Getz and Liu [GL] with no restriction on the global Schwartz functions.

\medskip\noindent
{\bf Acknowledgements}.
Confusions stemming from the shift $s+\frac 12$, which did not seem to agree with Ngo's shift [N3,BNS], did result in a number of communications with Erez Lapid for which I like to thank him. Similar gratitude are owed to Dihua Jiang, Wee Teck Gan, Shunsuke Yamana, David Kazhdan and Jayce Getz. I also like to thank Werner M\"uller, Sug Woo Shin and Nicolas Templier for their invitation to the Simons Symposium at Elmau, Germany, in April of 2016 and for the present proceedings. Parts of this paper were presented as a series of lectures at University of Minnesota where author was invited as an Ordway Distinguished Visitor during the Fall of 2016 and for which thanks are due to Dihua Jiang. Last but not least, I like to thank Wen--Wei Li for a numbers of helpful comments and communications after the first version of this manuscript was distributed, which in particular led to his insightful appendix [Li2] to this paper.

\bigskip\noindent
{\bf 1.\ $\gamma$--distributions.}

In this section we will introduce a distribution which will give local coefficients in the Langlands--Shahidi method, which are products of $\gamma$--factors, upon evaluation on the given representations. We note that this is parallel to what takes place with Godement--Jacquet $\gamma$--factors and its generalization by Braverman--Kazhdan.

We should point out that this section is rather independent of the rest of the paper. It is aimed at presenting a pair of distributions ($\gamma$--distributions), i.e., giving $\gamma$--factors upon evaluation on the corresponding representations, from two different methods, pointing out some similarities.

More precisely, let $F$ be a $p$--adic field of characteristic zero. Let $G$ be a connected reductive group over $F$. Fix a minimal parabolic subgroup $P_0$ over $F$ with a Levi decomposition $P_0=M_0N_0$, with $N_0$ the unipotent radical of $P_0$. We let $A_0$ be a maximal split torus of $G$, contained in $M_0$ as its split component. Let $P\supset P_0$ be a standard parabolic subgroup of $G$. Let $P=MN$ be a Levi decomposition with $M\supset M_0$ and $N_0\supset N$. Next, let $W=W(G, A_0)$. Finally, fix an irreducible admissible representation $\sigma$ of $M(F)$ and set
$$
I(\sigma)=\Ind^{G(F)}_{M(F)N(F)} \sigma\otimes \bold 1\tag1.1
$$
to denote the parabolically induced representation from $\sigma$. If $\Delta$ is the set of simple roots of $A_0$ in $N_0$, take $\theta\subset\Delta$ such that $M=M_\theta$. Let $\tw\in W$ be such that $\tw(\theta)\subset \Delta$. Fix a representative $w\in G(F)$ for $\tw$ which we will choose as in [Sh2,3]. We now define

$$
N_w:= N_0\cap wN^-w^{-1},\tag1.2
$$
where $N^-$ is the opposite of $N$. Given $f\in V(\sigma)$, the space of $I(\sigma)$, we have the intertwining operator
$$
A(\sigma,w)f(g):=\int\limits_{N_w(F)}\ f(w^{-1}ng)\ dn.\tag1.3
$$

Assume $P$ is maximal. Choosing representatives as in [Sh2,3], we let
$$
w=w_0=w_l\cdot w^{-1}_{l, M}. \tag1.4
$$
We now assume $w^{-1}_0Nw_0=N^-$, i.e., $P$ is self--associate. Then for each $n\in N(F)$, outside  a set of measure zero, there exist $m\in M(F)$, $n'\in N(F)$ and $\overline n\in\oN_{w_0}(F)=N^-(F)$ such that
$$
w^{-1}_0 n =mn'(\overline n)^{-1}.\tag1.5
$$
The intertwining operator (1.4) when evaluated at $g=\overline n_1\in\oN(F)$ can now be written as 
$$\aligned
A(\sigma, w_0)f(\overline n_1) &=\int\limits_{N(F)} \sigma(m)f((\overline n)^{-1}\overline n_1) dn\\
& \\
&=\Phi_\sigma(\overline n\longmapsto f((\overline n)^{-1}\overline n_1)),
\endaligned\tag1.6
$$
where $\Phi_\sigma$ is the measure defind via the bijection $n\mapsto\on$ of (1.5) by 
$$
\Phi_\sigma:=\sigma(m)dn.\tag1.7
$$
Thus
$$
A(\sigma, w_0)f=\Phi_\sigma * f.\tag1.8
$$

We now assume $G$ is quasisplit over $F$ in which case $P_0$ becomes a Borel subgroup of $G$ and $M_0=T$ is a maximal torus of $G$ with $T\supset A_0$. Let $\psi$ be a non-trivial additive character of $F$. Together with a fixed splitting $(G, B, T, \{x_\alpha\}_\alpha)$ of $G$, this defines a generic character of $U(F)$, still denoted by $\psi$. Assume $\sigma$ is $\psi$--generic. Let $\lambda$ and $\lambda'$ be the canonical Whittaker functional for $V(\sigma)$ and $V(w_0(\sigma))$, respectively. If $\lambda_M$ is a Whittaker functional for $\sigma$, then
$$
\lambda'(A(\sigma, w_0)f)=\int\limits_{\oN(F)}\lambda_M (A(\sigma,w_0)f(\overline n_1))\psi(\overline n_1)d\overline n_1,\tag1.9
$$
where $\psi(\overline{n_1}) : = \psi(w_0^{-1}\overline{n_1}w_0)$. The functional $\lambda$ is defined similarly. The local coefficient $C_\psi(\sigma)$ is defined by (cf. [Sh1,3])
$$
C_\psi(\sigma)^{-1}\lambda=\lambda'\cdot A(\sigma, w_0).\tag1.10
$$
We now formally define the $\psi$--Fourier transform of measure $\Phi_\sigma$ by
$$
\psi(\Phi_\sigma) :=\int\limits_{N(F)}\ \psi(\overline n)\sigma(m) d n.\tag1.11
$$
It is then the content of Proposition (9.4.15) of [Sh4] that
$$
\psi(\Phi_\sigma)\equiv C_\psi(\sigma)^{-1} \text{mod} (\ker(\lambda_M)).\tag1.12
$$
Moreover, using (1.5) for $u n u^{-1}$ for a fixed $u\in U_M(F)$, $U_M=U\cap M$, which reads
$$
w_0^{-1} unu^{-1}=w_0(u)mu^{-1}\cdot un'u^{-1}\cdot(u \overline n u^{-1})^{-1},\tag1.13
$$
$w_0(u)=w_0^{-1}uw_0$, one notices that
$$
\int\ \psi(\overline n) d n=\int\ \psi(u\overline n u^{-1}) d(unu^{-1})\tag1.14
$$
and thus the distribution
$$
\Phi_\psi=\psi(\overline n)d n\tag1.15
$$
is $\text{Int}(u)$--invariant (cf. [Sh4]).  We can therefore define an $\text{Int}(U_M)$--invariant distribution (1.15) which gives $C_\psi(\sigma)^{-1}$, a product of $\gamma$--factors within Langlands--Shahidi method, by evaluation at $\sigma$, i.e., 
$$\aligned
\Phi_\psi(\sigma):&=\int\limits_{N(F)}\psi(\overline n)\sigma(m) dn\\
&\equiv C_\psi(\sigma)^{-1}\mod\!\!(\ker(\lambda_M))
\endaligned\tag1.16
$$
by (1.12). We collect this as:

\medskip\noindent
{\bf(1.17)}\ {\bf Proposition}.
{\it  The value of the $\text{Int}(U_M(F))$--invariant distribution $\Phi_\psi=\psi(\overline n) d n$ on each $\psi$--generic irreducible admissible  representation $\sigma$ of $M(F)$ equals, up to $\ker(\lambda_M)$, with $C_\psi(\sigma)^{-1}$, a product of $\gamma$--factors of Langlands--Shahidi type.}

We will denote this distribution by $\Phi^{LS}_\psi$ to distinguish it from the next one, that of Godement--Jacquet (cf. [GJ,BK1]). (See the brief discussion at the beginning of Section 9 of [Sh4] for a quick review.) Recall that this involves the group $GL_n$ over local and global fields. Again $\gamma$--factors  are defined by means of a distribution
$$
\Phi_\psi(g)=\psi(tr(g)) |\det g|^n\ |dg|\tag1.18
$$
on $GL_n(F)$ as in Section 1.2 of Braverman--Kazhdan [BK1], where $F$ is a $p$--adic field. The $\gamma$--factor for the principal $L$--function of an irreducible admissible representation $\pi$ of $GL_n(F)$ is then equal to 
$$
\Phi_\psi(\pi)
=\int\limits_{GL_n(F)}\Phi_\psi(g)\pi(g). \tag1.19
$$
The distribution $\Phi_\psi$ is in fact $\text{Int}(G)$--invariant, $G=GL_n(F)$. We will denote this distribution by $\Phi^{GJ}_\psi$. One expects that the Godement--Jacquet theory will generalize to any reductive group and any representation of its $L$--group. Initial steps toward this are due to Braverman--Kazhdan [BK1,2], followed by Ngo [N1,2]. We summarize the discussion in this section as follows:

\bigskip
\noindent
{(1.20) \it
There are two distributions $\Phi^{LS}_\psi$ and $\Phi^{GJ}_\psi$, $\text{Int}(U_M(F))$-- and $\text{Int}(GL_n(F))$--invariant, respectively, whose values at an irreducible admissible $\psi$--generic representation $\sigma$ of $M(F)$ and an irreducible admissible representation $\pi$ of $GL_n(F)$, respectively, are a product of $\gamma$--factor for $\sigma$ of Langlands--Shahidi type (the inverse of the local coefficient) and the principal $\gamma$--factor for $\pi$, respectively.}

As mentioned above, one hopes that the approach of Godement--Jacquet can be generalized to an arbitrary reductive group and any finite dimensional representation of its $L$--group. We will make a quick survey of initial steps of this due to Braverman--Kazhdan [BK1,2] and Ngo[N1,2] in the next section.

\bigskip\noindent
{\bf 2.\ Monoids}

The purpose of this section is to provide the background on construction of monoids needed in Sections 3 and 6, as well as a brief comparison of the monoid structures with corresponding Levi subgroups appearing in Langlands--Shahidi method and its conjectural generalizations to Kac--Moody groups. 

As we pointed out in the Introduction and previous section, there has been initial efforts in generalizing the work of Godement--Jacquct [GJ] to all reductive groups and $L$--functions attached to arbitrary finite dimensional irreducible representations of their $L$--groups by Braverman--Kazhdan [BK1] and Ngo [N1,2]. The first step in this direction is a generalization of the space $M_n(F)$ of $n\times n$ matrices with entries in $F$ in which $GL_n (F)$ is the group of units. This generalization is called a 'monoid' which in characteristic zero was studied and classified by Vinberg [V], as we shall now explain. We refer to [N1,2,V] for other references and Section 9 of [Sh4] for a quick survey.

Let $k$ be an algebraically closed field of characteristic zero. Let $M$ be an irreducible affine algebraic normal variety over $k$ with an associative multiplication, i.e., a morphism
$$
\mu: M\times M\longrightarrow M
$$
of algebraic varieties so that $M$ is a semigroup. The null element $0$ or identity $1$ may or may not belong to $M$.
 If $1\in M$, then $M$ is called a {\bf{monoid}}. We let $G=G(M)$ denote the group of units of $M$ for $1$. The monoid $M$ is called {\bf{reductive}} if $G(M)$ is. The group $G=G(M)$ is never semisimple unless $M$ is a (semisimple) group. Let $G'=G_{der}$ be the derived group of $G=G(M)$. Then $G'\times G'$ acts on $M$ by
$$
(g_1, g_2)\cdot m=g_1 m g_2^{-1}.
$$
Let
$$
A:= M//G'\times G'\tag2.1
$$
be the $GIT$, geometric invariant theoretic quotient, of $M$ by $G'\times G'$. Then
$$
k[A]=k[M]^{G'\times G'}\hookrightarrow k[M].\tag2.2
$$
We will call $A$ the {\bf{abelianization}} of $M$. Let $\pi:M\longrightarrow A$ be the natural projection. We note that $k[A]\hookrightarrow k[M]$ is dual to $\pi$. We will assume $\pi$ is {\bf{flat}}, i.e., $k[M]$ is a free $k[A]-$module.

We now assume $G'$ is simply connected and let $T'$ be a maximal torus of $G'$ and $Z'$ be the center of $G'$.
Set
$$
T^+:=T'/Z'=T^{ad}\tag2.3
$$
the maximal torus of the adjoint group $G'/Z'$. We now set
$$
G^+:=T'\times G'/Z'.\tag2.4
$$
Next, fix a Borel subgroup $B'$ containing $T'$ and let $\{\alpha_1, \dots,\alpha_r\}$ be the set of simple roots of $(B', T')$. Then $T^+=T^{ad}$ can be identified with $\bG^r_m$ through the well defined map
$$
t\longrightarrow (\alpha_1(t), \dots, \alpha_r(t)),\tag2.5
$$
$t\in T^{ad}$, where $r$ is the (semisimple) rank of $G'$.

Now, let $\{\omega_1, \dots, \omega_r\}$ be the set of fundamental weights of $G'$, i.e. those dual to coroots
$$
\alpha^\vee_j=2\alpha_j / \kappa(\alpha_j,\alpha_j)\tag2.6
$$
through
$$
\kappa(\omega_i, \alpha^\vee_j)=\delta_{ij},\tag2.7
$$
when $\kappa$ is the Killing form and $\delta_{ij}$ the Kronecker delta function.

If $\rho_i$ is the fundamental representation attached to $\omega_i$ on the space $V_i$, we define an extention $\rho^+_i$ to $G^+$ by
$$
\rho^+_i(t,g)=\omega_i(w_0(t^{-1}))\rho_i(g),\tag2.8
$$
where $w_0$ is the long element of the Weyl group $W(G', T')$, $t\in T'$, $g\in G'$. We also extend $\alpha_i$ to $G^+$ by
$$
\alpha^+_i(t,g)=\alpha_i(t),\tag2.9
$$
$t\in T^{ad}$. We then get an embedding
$$
(\alpha^+, \rho^+): G^+\hookrightarrow \bG^r_m\times \prod^r_{i=1}GL(V_i).\tag2.10
$$
\medskip\noindent
{\bf(2.11)}\ {\bf Definition}.
{\it The closure of $G^+$ in}
$$
\bA^r\times\prod^r_{i=1}\End(V_i),
$$
$\bA=\bG_a$, {\it{is called the {\bf{Vinberg's universal monoid}}. It is denoted by $M^+$ and $G^+=G(M^+)$.
}}

We note that $M^+$ only depends on $G'$ as $G^+$ does.

Let $\pi^+:M^+\longrightarrow A^+$ be the abelianization of $M^+$. 
Vinberg's universal monoids theorem (Theorem 5 of [V]) determines every reductive monoid $M$ for which $G_{der}=G'$, $G=G(M)$, by means of its abelianization $A$ from the universal monoid $M^+$ which shares the same derived group $G'$ as $M$.

More precisely, assume there exists a morphism
$$
\varphi_{ab}: A\longrightarrow A^+.\tag2.12
$$
Then Vinberg gets $M$ as the {\bf{fibered product}} of $A$ and $M^+$ over $A^+$ by means of $\varphi_{ab}$, i.e.,
$$
M=A\times_{A^+}M^+:=\{(a,m^+)\in A\times M^+|\varphi_{ab}(a)=\pi^+(m^+)\}.\tag2.13
$$
We then have the commuting diagram
$$\CD
M @>{\varphi}>> M^+\\
{\pi}@VVV    @VV{\pi^+}V\\
A @>{\varphi_{ab}}>> A^+
\endCD\tag2.14
$$
where $\varphi$ is the projection of $M$ into $M^+$ in (2.13). In particular, $M$ is uniquely determined by $\varphi_{ab}$ if $0\in M$ (cf. [V]). 

We now consider the special case where $G'\setminus G\simeq\bG_m$. Let $\lambda:\bG_m\longrightarrow T^+\simeq\bG^r_m$ be a cocharacter of $T^+=T^{ad}$. Assume $\lambda$ is dominant, i.e., $\kappa(\lambda, \alpha^\vee)\ge 0$ for all simple roots $\alpha$. Then $\lambda$ extends to a map
$$
\varphi_{ab}:\bA^1\longrightarrow \bA^r\tag2.15
$$
and consequently by the main theorem defines a monoid $M^\lambda$. As we explain, $M^\lambda$ is the replacement of $M_n$ from the case of Godement--Jacquet [GJ], attached to the representation of the $L$--group of $G^\lambda:=G(M^\lambda)$ defined by the highest weight $\lambda$, a character of the dual torus $(T^{ad})^\wedge=\hat T^{sc}$. 

In fact, the dominant cocharacter $\lambda$ by duality defines a character
$$
\lambda:\hat{T}^{sc}\longrightarrow \bG_m\tag2.16
$$
giving the highest weight for an irreducible representation
$$
\rho_\lambda:\hat G^{sc}\longrightarrow GL(V_\lambda).\tag2.17
$$
If $\hat Z$ is the center of $\hat G^{sc}$, then $\hat Z$ acts on $V_\lambda$ by a character $\omega_\lambda$. The group $\hat G^{sc}$ is the covering group of the derived group of $(G^\lambda)^\wedge$, the dual group of $G^\lambda$, with the covering map
$$
F:\hat G^{sc}\longrightarrow (G^\lambda)^\wedge_{\der}.\tag2.18
$$
Then $F|\hat Z=\omega_\lambda$ which gives $\omega_\lambda(\hat Z)$ as the center of $(G^\lambda)^\wedge_{der}$. One then gets
$$
(G^\lambda)^\wedge=\bG_m\times F(\hat G^{sc})/\omega_\lambda(\hat Z)\tag2.19
$$
as detailed in [N1,2,Sh4]. In particular, one can directly compute $G^\lambda$ without calculating $M^\lambda$. Finally, we extend $\rho_\lambda$ from $\hat G_{sc}$ to $\rho^+_\lambda$ a representation of $(G^\lambda)^\wedge$ on $GL(V_\lambda)$ as in [N1,2].

\bigskip
\noindent{\bf{(2.20)  Examples}}. We refer to Section 9.3 of [Sh4] for a number of examples, including $G^\lambda$ for when $\lambda$ is the highest weight of a symmetric or exterior power representation of $GL_n(\bC)$, the $L$--group of $GL_n$. When $\lambda=\delta_2$ and $2\delta_1$, i.e., exterior or symmetric squares highest weights, which are among the cases appearing within Langlands--Shahidi method, the group $G^\lambda$ is precisely the Levi subgroup $L$ in a pair $(H, L)$, where $H$ is a simply connected semisimple group with a Levi subgroup $L$ giving the corresponding $L$--functions. We again refer to Section 9.3 of [Sh4] for other examples and the conjectural extension of Langlands--Shahidi method to infinite dimensional groups which we elaborate a bit here for the sake of completeness.

\bigskip
\noindent{\bf{(2.21)  Generalization to Kac--Moody groups}}. With notation as before, let $\lambda$ be a dominant cocharacter of $T^{ad}$. Denote by $G^\lambda$ the group of units of monoid $M^\lambda$ attached to $\lambda$ by Vinberg's universal monoids theory. Then as proved in Proposition 9.3.12 of  [Sh4], we can choose a complex adjoint Kac--Moody group $\tH$ and a maximal parabolic subgroup $\tP=\tG\tN$ of $\tH$ with a Levi subgroup $\tG$ such that the adjoint action $r$ of $\tG$ on $Lie(\tN)$ decomposes as $r=\oplus_i r_i$ with $r_1\cdot \eta$ containing $\rho_\lambda$, where $\eta:\hat G^{sc}\longrightarrow(\tG)_{der}$ is the covering map, and $\rho_\lambda$ is the representation of $\hat G^{sc}$ with highest weight $\lambda$ considered as a character of $\hat T^{sc}$ as in (2.17). The content of Proposition 9.3.11 of [Sh4] is then:
\bigskip
\noindent{\bf{Proposition}}.
{\it $\tG\simeq(G^\lambda)^\wedge$, where $(G^\lambda)^\wedge$ is the dual group of $G^\lambda$, the units of the monoid $M^\lambda$ attached to $\lambda$.}

\bigskip
\noindent{\bf{Remark}}.
We like to remark that at least in the  spherical case, the work of Patnaik [P] on unramified Whittaker functions when combined with  [BGKP] on Gindikin--Karpelevich formula, may lead to a definition of local coefficients in the case of loop groups over function fields (communications with Manish Patnaik). Whether any generalization of this to other Kac--Moody groups, non-spherical cases, or over number fields, i.e., a theory of local coefficients is possible, remains to be seen. It is definitely a very intriguing possibility. 

\bigskip
\noindent{\bf{(2.22)}  \bf{Symmetric powers for}} $\bold {GL_2}$. One important case where the monoid is easy to determine is that of symmetric power representations of $GL_2(\bC)$, the $L$--group of $GL_2$. Thus for any positive integer $n$, let $\lambda=n\delta_1$, the highest weight of $\Sym^n$, a representation of $SL_2(\bC)$. 
This is a character of maximal torus $GL_1(\bC)=\bG_m$ of $SL_2(\bC)$. We will then use $n\delta_1:\bG_m\longrightarrow\bG_m, n\delta_1(a)=a^n$, to define the corresponding cocharacter which is dominant, extending to a map $\varphi_{ab}: A^{n\delta_1}\longrightarrow A^+$, with $A^{n\delta_1}$ and $A^+$ abelianizations of $M^{n\delta_1}$ and $M^+$, respectively. 
Then as explained in [Sh4] both $A^+$ and $A^{n\delta_1}$ are isomorphic to $\bA^1=\bG_a$, $M^+=\End(\bA^2)$ and $\pi^+:M^+\longrightarrow A^+$ is simply the determinant map. It then follows from definition (2.13) that
$$M^{n\delta_1}=\{(a,m^+)\in\bA^1\times\End(\bA^2)| a^n=\det(m^+)\}$$
and
$$G^{n\delta_1}=\{(a,g)\in\bG_m\times GL_2| a^n=\det g\}.$$
It is proved in Lemma 9.3.21 of [Sh4] that $G^{n\delta_1}=GL_1\times SL_2$ if $n$ is even and $G^{n\delta_1}=GL_2$, otherwise.
The group of units $G^{n\delta_1}$ can easily be calculated from formula (2.19) directly. Moreover, when $n=2,3$, $G^{n\delta_1}$ is exactly the Levi subgroups in $Sp_4$ and $G_2$, giving symmetric square and cube $L$--functions for $GL_2$ within Langlands--Shahidi method (cf. [Sh1,3]). More precisely, the pairs are $(Sp_4, GL_1\times SL_2)$ and $(G_2, GL_2)$ for $n=2$ and 3, respectively. Note that both $Sp_4$ and $G_2$ are simply connected and thus in agreement with proposition 2.21.

\bigskip\noindent
{\bf 3.\ {A generalized Fourier/Hankel transform (d'apr\`es B.C. Ngo [N3])}}

Any generalization of Godement--Jacquet's work for $GL(n)$ to an arbitrary group $G$ and a representation of its $L$--group defined by a highest weight $\lambda$, beside the monoid $M^\lambda$ and its group of units $G^\lambda$, requires a space of Schwartz functions on $G^\lambda(F)$ which can be obtained by restricting smooth functions of compact support on $M^\lambda(F)$ to $G^\lambda(F)$. One crucial difference with the Godement--Jacquet theory is that in no other case $M^\lambda$ is smooth. Besides, if one can somehow define a space of Schwartz functions on $G^\lambda(F)$, one still needs an appropriate Fourier transform, and for global reasons, a Poisson summation formula for this transform over number fields.

To this end Ngo [N3] has now formulated a conjectural procedure to get this Fourier transform by defining it as a Hankel transform upon restriction to $F$--points of all the maximal tori of $G^\lambda$. The classical Hankel transform is the Fourier transform of radially symmetric functions on $\bR^2$, i.e., those which do not depend on $\theta$ in the polar coordinates, and thus defined as an integral transform whose kernel is a Bessel function. Fixing $\lambda$, let $\rho=\rho_\lambda$ be the representation of $\hat G^{\sc}$ defined by $\lambda$ as its highest weight via (2.16) and (2.17), and denote by $\rho_\lambda^+$, its extension to $(G^\lambda)^\wedge$ as in Section 2; also see Section 9.2 of [Sh4].

In [N3], Ngo explains how one can define a $\rho$--Fourier transform (Hankel transform) for a space of $\rho$--Schwartz functions on $T(F)$, where $T$ is a split torus of $G^\lambda$, using the projection $p_\rho$ of toric varieties 
$$
p_\rho: \bA^r\longrightarrow M_{T,\rho},\tag3.1
$$
and the standard Fourier transform $\sF$ for Schwartz functions $\sS(F^r)$ on $F^r$.

More precisely, following Ngo [N3], consider $\rho$ as a representation of $\hat T$ by restriction and let $\mu_1,\dots,\mu_r$ be its weights. The toric variety $M_{T,\rho}$ is the monoid defined by $\rho |\hat T$ or by the strictly convex cone generated by $\mu_1,\dots,\mu_r$ in $\Hom(\bG_m, T)\otimes\bR$.
The map $p_\rho$ is defined by
$$
p_\rho(x_1,\dots,x_r)=\mu_1(x_1)\dots\mu_r(x_r).\tag3.2
$$
Note that $T\subset M_{T,\rho}$ and $\bG_m^r\subset \bA^r$ are open embeddings giving units of $M_{T,\rho}$ and $\bA^r$. Consequently, $p_\rho$ gives a map from $\bG^r_m\rightarrow T$. Let $S$ be its kernel. Then $S(F)$ acts on $\sS(F^r)$ by
$$
(s\cdot f)(x)=f(s^{-1}\cdot x),\tag3.3
$$
where $s\in S(F)$, $x\in F^r$ and $f\in\sS(F^r)$, with multiplication of $\bG^r_m$ coordinate--wise on $\bA^r$.

The $\rho$--Schwartz space $\sS_\rho(T(F))$ of $T(F)$ is just the coinvariants of $\sS(F^r)$ under the action of $S(F)$.
The push forward
$$p_{\rho,!}:\sS(F^r)\longrightarrow\sS_\rho(T(F))\tag3.4
$$
is then simply defined by integration over the fibers, i.e., 
$$\aligned
p_{\rho,!}(f)(x)&=\int\limits_{p^{-1}_\rho(x)} f(y)dy\\
\\
&=\int\limits_{S(F)} f(y_0 s)d s,
\endaligned\tag3.5
$$
with $x=p_\rho(y_0)$ and, where $ds$ is a Haar measure on $S(F)$. 

The $\rho$--Fourier transform (Hankel transform)
$$\sF_\rho:\sS_\rho(T(F))\longrightarrow \sS_\rho(T(F))\tag3.6
$$ 
is now defined formally to satisfy
$$
\sF_\rho(p_{\rho,!} f)=p_{\rho,!}(\sF f),\tag3.7
$$
$f\in\sS(F^r)$, where $\sF$ is the euclidean (standard) Fourier transform on $\sS(F^r)$, i.e., the following diagram commutes:
$$\CD
\sS(F^r) @>{\sF}>> \sS(F^r)\\
{p_{\rho, !}}@VVV    @VV{p_{\rho, !}}V\\
\sS_\rho(T(F)) @>{\sF_\rho}>> \sS_\rho(T(F))
\endCD\tag3.8
$$

In [N3], Ngo denotes $\sF_\rho$ by $\sJ_\rho$, the Hankel transform, which is uniquely defined by (3.7). It can formally be given by
$$
\sJ_\rho(\theta)=J_{T, \rho}\ast {\check{\theta}},\tag3.9
$$
$\theta\in\sS_\rho(T(F))$, ${\check{\theta}}(x)=\theta(x^{-1})$, in which 
$$
J_{T,\rho}=p_{\rho,!}J_{\std},\tag3.10
$$
where the kernel $J_{\std}$ of euclidean (standard) Fourier transform is defined by
$$
J_{\std}(x)=\psi(x)|x_1|\dots|x_r|,\tag3.11
$$
with $\psi$ a non--trivial additive character of $F$, giving one of $F^r$ by $\psi\cdot tr$, and defining the standard Fourier transform by
$$
\sF(\phi)=J_{\std} * {\check{\phi}},\tag3.12
$$
$\phi\in\sS(F^r)$.

Ngo [N3] then extends this to non--split tori, first treating the case where restriction of $\rho$ to $\hat T$ is multiplicity free and then the general case. In fact, the situation is more complicated, but becomes more manageable if the weights $\mu_i$ of $\rho |\hat T$ are multiplicity free. The action of the Galois group $\Gamma$ on weights can then be expressed through a homomorphism
$$\rho_{_{\Gamma}} :\Gamma\longrightarrow\frak S_r\tag3.13
$$
such that $\sigma(\mu_i)=\mu_{\rho_{_{\Gamma}}}(\sigma)$, $\sigma\in\Gamma$, where $\frak S_r$ is the symmetric group in $r$--letters, giving the action of $\sigma$ on weights of $\rho |\hat T$. One can now put a $\Gamma$--structure on $\bA^r$ such that the map
$$p_\rho :\bA^r\longrightarrow M_{T,\rho}\tag3.14
$$
defined by (3.2), becomes $\Gamma$--equivalent. The construction in the split case now carries over, leading to construction of $\sS_\rho(T(F))$ and $\sF_\rho$ in this case.

For the case where there are multiplicities, we refer to Ngo's construction in second half of Section 5.2 in [N3].

As we discussed earlier, we are interested in the case that the torus $D=G/G'$, $G'=[G, G]$, is one--dimensional and, $D=\bG_m$. Thus
$$
0 @>>> G' @>>> G @>c >> \bG_m @>>> 0.\tag3.15
$$
We then have
$$
0 @>>> \Lambda' @>>> \Lambda @>c >> \Lambda_D @>>> 0,\tag3.16
$$
where $\Lambda=\Hom(\bG_m, T)$, $\Lambda'=\Hom(\bG_m, T')$ and $\Lambda_D=\Hom(\bG_m, \bG_m)\simeq\bZ$. We set $\Lambda^+\subset\Lambda$ to be the dominant cocharacters, i.e., those for which $c(\Lambda^+)\simeq \bZ^+$, non--negative integers. Choose $\lambda\in\Lambda^+$, projecting to $1\in\bZ^+$. Let $\rho$ be the irreducible representation
$$
\rho:\hat G\longrightarrow GL(V\rho)\tag3.17
$$
defined by the highest weight $\lambda$ as before, where $V_\rho$ denotes the space of $\rho$. Here we use $\hat G:=(G^\lambda)^\wedge$. 
More precisely, we note that here $\lambda$ is in fact the highest weight of the representation $\rho_\lambda$ of the simply connected group $\hat G^{\sc}$ as in (2.17) which is lifted to one of $(G^\lambda)^\wedge=\hat G$ defined in (2.19). The action of $\hat c(t)$ of the center $\bC^*$ of $(G^\lambda)^\wedge$ on $V_\rho$ is scalar multiplication by $t$, where
$$\hat c:\bC^*\longrightarrow (G^\lambda)^\wedge=\hat G\tag3.18
$$
is the dual map, embedding $\bC^*$ as the center of $(G^\lambda)^\wedge=\hat G$. This can simply be stated as
$$\rho\cdot\hat c(t)=t.\tag3.19
$$

For our purposes we like to recall basic functions. Let $\pi$ be an irreducible unramified representation of $G(F)$, i.e., $\pi^{G(O_F)}\neq \{0\}$. Then $\pi$ is parameterized by an element $\sigma=\sigma_\pi\in\hat G$. More precisely, $\pi$ will be uniquely determined by the conjugacy class of $\sigma$ and $\sigma$ is semisimple. The Langlands $L$--function attached to $\pi$ and $\rho$ is defined by 
$$\aligned
L(s,\pi,\rho) &=\det(1-\rho(\sigma)q^{-s})^{-1}\\
& =\sum^\infty_{d=0} tr(\sym^d\rho(\sigma))q^{-ds},\endaligned\tag3.20
$$
$s\in \bC$, a formal sum which converges for $Re(s)\gg 0$.

Let $\sH=\sH(G(F),K)$, $K=G(O_F)$, be the spherical Hecke algebra defined by $K$. It is the space of bi--$K$--invariant complex functions of compact support on $G(F)$ under convolution.
Satake isomorphism
$$\Sat: \sH(G(F), K)\simeq\bC[\hat G]^{ad(\hat G)}\tag3.21
$$
is a canonical one between $\sH(G(F),K)$ and the space of regular functions on $\hat G$, invariant under $\hat G$--conjugation. Then
$$tr_\pi(\phi)=(\Sat (\phi))(\sigma_\pi),\tag3.22
$$
where
$$tr_\pi(\phi)=tr(\pi(\phi))\tag3.23
$$
in which
$$\pi(\phi)=\int \, \phi(g)\pi(g) dg\tag3.24
$$
is an operator of finite rank.

If $\eta_G$ is half the sum of roots in $B$, i.e.,
$$\eta_G=\frac 12\sum_{\alpha>0}\alpha,\tag3.25
$$
we can define
$$\bL^\rho(s)=\sum^\infty_{d=0}\, \phi^{\sym^d(\rho)} q^{-ds}.\tag3.26
$$
Then
$$\aligned
L(s,\pi,\rho)&=\sum^\infty_{d=0}\,tr_\pi(\phi^{\sym^d(\rho)})q^{-ds}\\
\\
&=tr_\pi(\bL^\rho(s))\\
\\
&=tr(\bL^\rho, \pi\otimes |c|^{s+\langle \eta_G, \lambda\rangle}),\endaligned\tag3.27
$$
where
$$\bL^\rho : =\bL^\rho (-\langle \eta_G, \lambda\rangle ),\tag3.28
$$
the {\bf basic function} attached to $\rho$.

The traditional Satake isomophism for a function $\phi\in\sH(G(F),K)$ is given by its constant term $\phi_N$, namely
$$\phi_N(t)=\delta_B(t)^{1\!/\! 2}\int\limits_{N(F)}\phi(tn)dn,\tag3.29
$$
where $N$ is the unipotent radical of the Borel subgroup $B=TN$.

When the function $\phi=\phi^\lambda$ corresponds to the trace of a representation $\rho$ of highest weight $\lambda$ as a Schur spherical function, then
$$\phi^\lambda_N\, (\mu)=m_\mu (V_\rho),\tag3.30
$$
where for each weight $\mu$, $m_\mu(V_\rho)$ gives its multiplicity in $V_\rho$. In particular, for our basic function $\bL^\rho$
$$\bL^\rho_N (\mu)=m_\mu (\sym(\rho))\tag3.31
$$
in which $m_\mu(\sym(\rho))=m_\mu(\sym^d(\rho))$, where $d=c(\mu)$, $c:\Lambda\rightarrow \bZ$ as before.

Let $\nu_1,\dots,\nu_r$ be a basis of $V_\rho$, given by weights $\mu_1,\dots,\mu_r$ of $\rho$. If $\mu\in\Lambda$ has degree $d=c(\mu)$, then the multiplicity $m_\mu(\sym^d(V_\rho))$ can be given by a partition function attached to certain multisets on $\{\nu_1,\dots,\nu_r\}$. More precisely, let $p^{\mu_1,\dots,\mu_r}(\mu)$ be the number of ways one can find multisets $\{\nu_{i_1},\dots,\nu_{i_d}\}$, $d=c(\mu)$, for which corresponding weights satisfy the partition
$$\mu_{i_1}+\cdots+\mu_{i_d}=\mu\tag3.32
$$
of $\mu$. We recall that in a multiset multiplicities are allowed.

If ${\bold 1}_{O^r}$ is the characteristic function $O^r$, $O=O_F$, then it can be shown that fixing measures on $(F)^r$ and $T(F)$ suitably, its push--forward $p_{\rho, !}({\bold 1}_{O^r})$ equals to this partition function and thus 
$$p_{\rho,!}(\bold 1_{O^r})=p^{\mu_1,\dots,\mu_r}=\phi^\lambda_N,\tag3.33
$$
where $p_\rho$ is as in (3.2).

Applying (3.33) to a given weight $\mu$ will now give the weight multiplicity of $\mu$ in each $\sym^d(\rho)$, $d=c(\mu)$ which vanishes if $d\ne c(\mu)$. In particular, the push--forward  $p_{\rho, !}(\bold 1_{O^r})$ equals to $\bL^\rho(0)$ evaluated on $T(F)$. But
$$\bL^\rho(0)(g)=|c(g)|^{\langle \eta_G,\lambda\rangle}\bL^\rho(g),\tag3.34
$$
$g\in G(F)$. Thus
$$p_{\rho,!}({\bold 1}_{O^r})(t)=|c(t)|^{\langle \eta_G,\lambda\rangle}\bL^\rho(t),\tag3.35
$$
for all $t\in T(F)$. 

\bigskip\noindent
{\bf 4.\ {Fourier transforms of Braverman--Kazhdan}}

There is one case outside that of Godement--Jacquet [GJ] where the Fourier transform is explicitly given and that is the doubling method of Piatetski--Shapiro and Rallis [GPSR,PSR]. The idea was generalized and put in the context of this approach by Braverman and Kazhdan in [BK2] and was later elaborated on by Wen--Wei Li [Li]. As we will see later the Fourier transform in this case is given by a normalized intertwining operator [LR]. In this section we will review and reinterpret the normalizing factors in [BK2]. This will enable us to prove their equality with those of Piatetski--Shapiro and Rallis [GPSR,PSR], as well as Lapid and Rallis [LR], later in Section 6. These normalizing factors are not the standard ones suggested by Langlands [A,La,Sh3], but rather a singular version [LR] of local coefficients defined in [Sh1]. In what follows we mix the notation and results from [BK2] and [Li].

Let $F$ be a $p$--adic field with $O=O_F$ its ring of integers. Let $H$ be a split connected reductive group over $F$. We fix a parabolic subgroup $P$ with a fixed Levi decomposition $P=MU_P$, where $U_P$ is the unipotent radical of $P$. Set
$$M_{ab}\,:\,= M/M_{\der}\simeq P/P_{\der},\tag4.1$$
where $M_{\der}=[M,M]$ and $P_{\der}=[P,P]$. Let
$$X_P:\,= P_{\der} \setminus H.\tag4.2$$
It is a $M_{ab}\times H$--space by
$$(P_{\der} y)\cdot (\om,h)=P_{\der} m^{-1}y h.\tag4.3$$
Then
$$P\to P/P_{\der}\simeq M_{ab}$$
and $P\hookrightarrow H$ gives an embedding
$$P\hookrightarrow M_{ab}\times H$$ and one can identify $X_P$ with $P\setminus M_{ab}\times H$.
For $F$--points we have
$$\aligned
X_P(F)&=(P_{\der} \setminus H)(F)\\
&=P_{\der}(F)\setminus H(F),\endaligned$$
if $H^1(P_{\der})$ is trivial. This will be the case if $H_{\der}$ is simply connected (cf. next subsection); but it will also be the case for classical groups where $P$ is the Siegel parabolic for which $M=GL_n$, the case of interest to us is Sections 5 and 6. As is shown in [Li], it also follows if $M(F)$ projects onto $M_{ab}(F)$, which is clearly true if $H^1(P_{\der})$ is trivial.

Given a character $\chi$ of $M_{ab}(F)$ and a representation $\pi$ of $H(F)$, Frobenius reciprocity implies
$$\aligned
&\Hom_{(M_{ab}\times H)(F)}(\chi\otimes\pi, C^\infty(X_P))=\Hom_{\ P(F)}(\chi\otimes \pi, \delta^{1\!/\!2}_P)\\
\\
&=\Hom_{P(F)}(\pi, \tchi \otimes\delta^{1\!/\!2}_P)=\Hom_{H(F)}(\pi,\Ind^{H(F)}_{P(F)}\tchi)
\endaligned\tag4.4
$$
where
\medskip
\noindent
$\ds(4.5)\hskip1.5truein 
\tchi : M(F)\longrightarrow M_{ab}(F) \overset\chi\to  \longrightarrow \bC^{\ast}$.
\medskip
Now take $\pi$ to be the right action of $H(F)$ on $C^\infty_c(X_P)$. Then the map
$$
\tchi\otimes\xi \longmapsto \int\limits_{M_{ab}(F)} \chi(m)(m\xi)(\cdot) dm
\tag4.6
$$
which picks up the $\chi$--component of $C^\infty(X_P)$ gives an $H(F)$--map
\medskip
\noindent
$\ds(4.7)\hskip1.5truein C^\infty_c(X_P)\longrightarrow {\Ind}^{H(F)}_{P(F)} \tchi$.

\bigskip

\noindent
{\bf The Schwartz space of Braverman--Kazhdan.}

We now assume $H_{der}$ is simply connected and fix a non-trivial additive character $\psi$ of $F$.

We will first define a Schwartz space in this case which turns out to be the $\rho$--Schwatz space when $\rho$ is the standard representation of an appropriate $L$--group.

Let $P$ and $Q$ be two parabolic subgroups of $H$, sharing the same Levi subgroup $M$. Then in [BK2], Braverman and Kazhdan define an (intertwining) map
$$\sF_{Q|P}=\sF_{Q|P,\psi} : L^2(X_P)\longrightarrow L^2(X_Q)\tag4.8$$
which is $(M_{ab}\times H)(F)$--equivariant and an isometry. As we range the parabolic subgroup, the family of maps satisfy
$$\sF_{R|Q}\cdot \sF_{Q|P}=\sF_{R|P}\tag4.9$$
and $\sF_{P|P}={\text id}$, and thus
$$\sF_{P|Q}\cdot \sF_{Q|P}={\text id}.\tag4.10$$

Using the map (4.7), they therefore behave like normalized standard intertwining operators between induced representations.

The {\bf Schwartz space} $\sS(X_P)$ is defined as:
$$\sS(X_P): =\sum\limits_Q \sF_{P|Q} (C^\infty_c(X_Q)),\tag4.11$$
where $Q$ runs over parabolic subgroups sharing the same Levi as $P$.
The space is a smooth $(M_{ab}\times H)(F)$--representation through the action on different $X_Q$ of $(M_{ab}\times H)(F)$.

We now explain how $\sF$ is defined and how the diagram below commutes:
$$\CD
\sS(X_Q) @>{\sF _{P|Q}}>> \sS(X_P) \\
@VVV   @VVV \\ 
\Ind^{H(F)}_{Q(F)} \tchi @>> N_{P|Q}> \Ind^{H(F)}_{P(F)} \tchi  
\endCD\tag4.12
$$ 
In particular, $\sF_{P|Q}$ projects to a normalized standard intertwining operator $N_{P|Q}$. The vertical arrows are defined by (4.6). We remark that there are several ways to normalize the standard operator. This one is not the standard one used in the trace formula (cf. [A,L,Sh3]). We will explain this later.

The map $\sF$ is defined in two steps in [BK2] as follows:

\noindent
{\bf Step 1. The map} $\sR_{P|Q}$.

This projects to a non-normalized standard intertwining operator
\medskip
\noindent
$\ds (4.13)\hskip1.5truein
 J_{P|Q} : \Ind^{H(F)}_{Q(F)} \tchi\longrightarrow \Ind^{H(F)}_{P(F)} \tchi$.
\medskip
We recall that $P$ and $Q$ share the same Levi subgroup $M$.

Let $Z_{P,Q}\subset X_P\times X_Q$ be the image of $H$ under projection to $X_P$ and $X_Q$. Given $f\in C_c^\infty(X_Q(F))$, define
$$\sR_{P|Q}(f)(x)=\int\limits_{(x,y)\in Z_{P,Q}(F)}\, f(y)dy.\tag4.14$$
We now interpret this integral as follows. There exist elements $h\in H(F)$, $p\in P_{\der}(F)$ and $q\in Q_{\der}(F)$, such that $h=px=qy$. The elements $p$ and $q$ are unique up to an element in $(P_{\der}\cap Q_{\der})(F)$. We can therefore consider the integral of $f(px)=f(qy)$ over $(P_{\der}\cap Q_{\der})(F)\setminus P_{\der}(F)$, in which the measure $dp=dh/dx$, where $dh$ and $dx$ are invariant measures on $H(F)$ and $X_P(F)$, respectively. The function $\sR_{P|Q}(f)$ is now given by
$$\sR_{P|Q}(f)(x)=\int\limits_{(P_{\der}\cap Q_{\der})(F)\setminus P_{\der}(F)}\ f(px)dp.\tag 4.15$$
Note that $(P_{\der}\cap Q_{\der})(F)\setminus P_{\der}(F)\simeq (U_P\cap U_Q)(F)\setminus U_P(F)$ and therefore $\sR_{P|Q}$ simply projects to the standard intertwining operator:

\medskip
\noindent
$$J_{P|Q} : {\Ind^{H(F)}_{Q(F)}}{(\tchi)}\longrightarrow  {\Ind^{H(F)}_{P(F)}}{(\tchi)}$$
\medskip
\noindent
via
$$\CD
C^\infty_c(X_Q) @>{\sR _{P,Q}}>> C^\infty(X_P) \\
@VVV   @VVV \\ 
\Ind^{H(F)}_{Q(F)} \tchi @> J_{P|Q}>> \Ind^{H(F)}_{P(F)} \tchi  
\endCD\tag4.16
$$ 
for each $\tchi$--component of corresponding Schwartz spaces $\sS(X_Q)$ and $\sS(X_P)$.

{\bf Step 2. The normalizing factor.} 

Braverman and Kazhdan [BK2] define these factors on $C_c^\infty(X_Q(F))$ as a distribution, which when projected via diagrams (4.12) and (4.16) becomes a normalizing factor for $J_{P|Q}$. We shall now explain.

Let $T$ be a split torus in $H$. We have the lattice of cocharacters
$$\aligned
\Lambda_*(T)&=\Hom (\bG_m, T)\\
&=\Hom(\hat{T},\bG_m),\\
&=\Lambda^*(\hat T),
\endaligned\tag4.17
$$
where $\hat{T}$ is the dual group of $T$. For us eventually $\hT=\hM_{ab}=Z(\hM)$ with obvious notation for the objects.

Let $L=\bigoplus\limits^k_{i=1}L_i$ be a graded finite dimensional representation of $\hT$ defined by a collection of elements $\lambda_1,\dots,\lambda_k\in X^*(\hT)=X_*(T)$ identified as $\hT$--eigenbasis for $L$. We assign integers $n_1,\dots,n_k$ to these eigenspaces which we will be more specific in the case of interest to us later in Section 6. We allow multiplicity among $\lambda_i$ and $n_i$. We let $s_i=n_i/2$.

Considering $\lambda_i:\bG_m\to T$, we can consider its push--forward $(\lambda_i)!$ on the space of distributions on $\bG_m$. Thus
$$(\lambda_i!)(\eta)(\varphi)=\eta(\varphi\cdot\lambda_i),\tag4.18$$
for each distribution $\eta$ on $F^*$ and each function $\varphi$ on $T(F)$ whenever it makes sense; thus $\lambda_i!$ transfers distributions on $F^*$ to those on $T(F)$.

We will now restrict ourselves to a specific distribution on $F^*$. Let $\psi$ be a non-trivial additive character of $F$. Fix a self--dual measure $dx$ with respect to $\psi$. Let $s\in \bC$. Define the distribution
$$\eta=\eta^s_\psi : =\psi(x)|x|^s|dx|.\tag4.19$$
The distribution $\eta^s_\psi=\eta$ can be integrated  against a character $\chi\in \hF^*$, i.e.,
$$\langle\eta, \chi\rangle =\int \chi(x)\psi(x) |x|^s |dx|.\tag4.20$$
The reader who is familiar with Gauss sums realizies that $\langle\eta, \chi\rangle$ converges as a principal value integral and defines a rational function via the Laurent power series:
$$
M_{\eta,\chi}(z)=\sum\limits^\infty_{n=-\infty} z^n \int\limits_{|x|=q^{-n}} \chi(x)\psi(x) |dx|,\tag4.21
$$
where $z=q^{-s}$. Setting $z=1$ in 
$$M_{\eta,\chi}(z)=\sum\limits^\infty_{n=-\infty} z^n \int\limits_{|x|=q^{-n}}\eta(x)\chi(x),$$
 we get a rational function $M(\eta)(\chi)=M_{\eta,\chi}(1)$. By formula (4.20), it is clear that
$$M(\eta^s_\psi)(\chi)=\gamma(s,\chi,\psi),\tag4.22$$
the Hecke--Tate $\gamma$--function. The distribution $\eta^s_\psi$ which makes (4.21) convergent is an example of what Braverman--Kazhdan call a ''good distribution''. We refer to [BK2] for a discussion of these distributions and their properties.

Given our graded $\hT$--representation $L$, one can now define as in [BK2] a distribution
$$\aligned
\eta_{_{L,\psi}}\ :\  &=\eta^{s_1, s_2,\dots,s_k}_{\lambda_1,\lambda_2,\dots,\lambda_k,\psi}\\
&=(\lambda_1!)(\eta^{s_1}_\psi)\ast(\lambda_2!)(\eta^{s_2}_\psi)\ast\cdots\ast(\lambda_k!)(\eta^{s_k}_\psi),
\endaligned\tag4.23
$$
where $s_i={n_i\over 2}$. This is a distribution on $T(F)$.

It is shown in [BK2] that $\langle \eta_{_{L,\psi}},\chi\rangle$ is again given by a rational function $M(\eta_{_{L,\psi}})(\chi)$,
\vskip.01truein\noindent
where $\chi$ is a character of ${T(F)}$. Moreover,
$$M(\eta_{_{L,\psi}})=M((\lambda_1!)(\eta^{s_1}_\psi))\dots M((\lambda_k!)(\eta^{s_k}_\psi)).\tag4.24$$
Note that since $\chi\in\widehat{T(F)}$, $\chi.\lambda_i\in\hat{F}^*$ and thus
$$\aligned
M((\lambda_i!)(\eta^{s_i}_\psi))(\chi)&=M((\eta^{s_i}_\psi))(\chi.\lambda_i)\\
&=\gamma(s_i, \chi \cdot\lambda_i,\psi)
\endaligned\tag4.25
$$
and therefore
$$M(\eta_{_{L,\psi}})(\chi)=\prod^k_{i=1}\gamma(s_i, \chi\cdot\lambda_i,\psi).\tag4.26$$
We record this as

\bigskip
\noindent{\bf{(4.27) Proposition}}.
{\it Let $\eta_{_{L,\psi}}$ be as in (4.23). Let $M(\eta_{_{L,\psi}})$ be the rational function attached to it. Then
$$M(\eta_{_{L,\psi}})(\chi)=\prod^k_{i=1}\gamma(s_i, \chi\cdot\lambda_i,\psi),$$
where $\chi\in\hat{F}^*$.}

Now assume $T=M_{ab}$. One can then consider the convolution $\eta_{_{L,\psi}}\ast\varphi$ of the

\noindent
distribution $\eta_{_{L,\psi}}$ on $M_{ab}(F)$ with any function $\varphi\in C^\infty_c(X_P)$ on $M_{ab}(F)$ to get

\noindent
another distribution on $M_{ab}(F)$. Given $\chi\in \widehat{M_{ab}(F)}$, we calculate
$$
\aligned
\int\limits_{M_{ab}(F)}(\eta_{_{L,\psi}}\ast\varphi(\cdot))(m)\chi(m)dm
&=\int\limits_{(M_{ab}(F))^2}\eta(m_1)((m_1^{-1}m)\cdot\varphi)(\cdot)\chi(m)dm_1dm\\ \\
(4.28)\hskip1.75truein &=\int\limits_{M_{ab}(F)}\eta(m_1)\chi(m_1)dm_1\cdot \int\limits_{M_{ab}(F)} (m\cdot\varphi)(\cdot)\chi(m) dm\\ \\
&=M(\eta_{_{L,\psi}})(\chi)\int\limits_{M_{ab}(F)} (m\cdot\varphi)(\cdot)\chi(m) dm,
\endaligned
$$
where $\eta=\eta_{_{L,\psi}}$. We thus have

\bigskip
\noindent{\bf{(4.29) Proposition}}.
{\it The convolution by $\eta_{_{L,\psi}}$ of elements of $C^\infty_c(X_P)$ covers multiplication by the rational function $M(\eta_{_{L,\psi}})(\chi)$ on ${\Ind}^{H(F)}_{P(F)}\tchi$ for each $\chi\in \widehat{M_{ab}(F)}$, via (4.12).}

We will conclude the definition of Fourier transform $\sF_{P|Q}$ later in Section 6.
\bigskip\noindent
{\bf 5.\ {Connection with doubling method of Piatetski--Shapiro and Rallis}}

We shall now connect the Fourier transform of the previous section to the ''doubling method'' of Piatetski--Shapiro and Rallis [GPSR,PSR,LR] which is now completed for all classical groups in [L,Y1,2] including their inner forms in the unitary case;  and provides a theory of standard $L$--functions for these groups. (We refer to [Gan] for the case of double cover of symplectic groups.) In section 6, we determine the corresponding ''basic functions'' and show that they are fixed by the Fourier transform.

We will also address the issue of required shift introduced in [BNS,N3]. As we explain, a consistent comparison with the case of Godement--Jacquet [GJ] principal $L$--functions for $GL_n$, will only be possible if their $L$--function is also studied in the frame work of doubling method (cf. [Y2]). 

To explain doubling, we follow [LR]. Let $E$ be a local field ($p$--adic) of characteristic zero with an automorphism $\theta$ of order 1 or 2 with a fixed field $F$. Let $|\ |=| \ |_E$ be the normalized absolute value of $E$, given by $|\varpi_E|_E=q^{-1}_E$, where $q_E$ is the order of $O_E/P_E$.  Let $h$ be sesqui--linear $E$--valued (linear $F$--valued if $E=F$) form on an $n$--dimensional vector space $V$ over $E$ for which
$$\theta(h(v,u))=\varepsilon h(u,v)\quad (u,v\in V),\tag5.1$$
where $\varepsilon\theta(\varepsilon)=1$, $\varepsilon\in E^*$. If $h$ is non--degenerate, then it is either symmetric or skew--symmetric (symplectic) when $E=F$, $\varepsilon=1$ or $-1$, respectively; or $\varepsilon$--hermitian if $E\ne F$. By Hilbert Theorem 90, we may assume $\varepsilon=1$ in the hermitian case. We will also allow $h=0$.

We then consider the pair $(V,h)$ and let $G$ be the isometry group of this pair as an algebraic group over $F$, $G(F)\subset GL_n(E)$ or $G(F)\subset GL_n(F)$, depending whether the form $h$ is hermitian or not. When $h=0$, then $G(F)=GL_n(F)$. We refer to [GPSR,LR] for detail.

Doubling introduces a doubling of our data. We consider the space $W=V\oplus V$ together with the form $h^\square=h\oplus -h$, i.e.,
$$h^\square((v_1,v_2), (v'_1,v'_2))=h(v_1, v'_1)-h(v_2, v'_2).\tag5.2$$
Let $H=G^\square$ be isometry group of $h^\square$. We identify $G\times G$ with the subgroup of $H$, preserving $V_1=V\oplus \{0\}$ and $V_2=\{0\}\oplus V$. 
If $h\ne 0$, then
$$V^d=\{(v,v)| v\in V\},\tag5.3$$
is a maximal totally isotropic subspace of $(W,h^\square)$, since $h^\square((v_1,v_1),(v_2,v_2))=0$. Consequently, the stabilizer of $V^d$ will be a maximal parabolic subgroup $P$ whose Levi subgroup $M\simeq GL(V^d)=GL_n$. In particular, $H$ will be a quasisplit group over $F$.

As in previous section the main object through which Schwartz functions can be defined is an induced representation. More precisely, let $\chi$ be a character of $E^*$ and $s\in \bC$. Let $\chi_s=\chi\cdot |\  |^s_E$. Define a character of $M(F)=GL(V^d)(F)$ by $\chi_s\cdot \Delta$, where $\Delta=\det$ under the isomorphism $M\simeq GL(V^d)$. Set
$$I(s,\chi): = \Ind^{H(F)}_{P(F)}\chi_s\cdot\Delta,\tag5.4$$
the {\bf normalized} induced representation from $(\chi_s\cdot\Delta)\otimes {\bold 1}$, a representation of $M(F)U_P(F)$.

\medskip\noindent
(5.5)\ {\bf Remark}. 
We note that the modulus character of $P$ can be expressed as a power of $\det$. When $h$ is symplectic $\delta_P=\Delta^{n+1}$, while $\delta_P=\Delta^{n-1}$ for $h$ symmetric and $n$ even. These exponents are just the values $\langle 2\rho, \varepsilon_1\rangle$, where $\rho$ is half--sum of positive roots in $G$ and $\varepsilon_1$ is the first coordinate function for $T_H$, the standard diagonal torus of $H$.

The main goal of doubling method was to develop the theory of standard $L$--functions for classical groups along the lines of Godement--Jacquet [GJ], i.e., by means of matrix coefficients, and thus for any irreducible admissible representation $\pi$ of $G(F)$. As pointed out earlier, this is the $L$--function attached to the standard representation of $\hat G$. In the cases of classical groups the embedding $\hat G\hookrightarrow GL_N(\bC)$, with $N$ minimal, will give the standard representation of $\hat G$ upon restriction from that of $GL_N(\bC)$.

What was discussed in Section 4 now gives a new meaning to the introduction of $I(s,\chi)$ and the corresponding intertwining operator. Recall that $I(s,\chi)$ was originally introduced by Piatetski--Shapiro and Rallis in [GPSR], to prove the finiteness of poles of global standard $L$--functions through the theory of Eisenstein series.

While the analogy with the work of Godement--Jacquet [GJ] in the case of $GL_n$ was apparent, for example, through the use of matrix coefficients, its connection to the geometry involved and theory of reductive monoids discussed in previous sections was not transparent and clearly not addressed. We will now address this further, complementing the writings of Wen--Wei Li [Li] and Braverman--Kazhdan [BK2].

One of the main motivation for a local theory as proposed by Braverman--Kazhdan [BK1] and Ngo [N1,3], is to define the local $\gamma$--factors needed to derive the global functional equation. In the case of doubling these factors are now defined in every case, starting by Piatetski--Shapiro and Rallis's announcement [PSR], with complete proofs given by Lapid and Rallis in [LR], and extended, using the methods in [GPSR], to unitary groups by Jian--Shu Li [L] and further extended by Yamana in [Y1,2]. (The work in [GPSR] carries some typos which thankfully are corrected in [LR].) The main tool in definition is a normalized intertwining operator originating from $I(s,\chi)$ which we shall now explain:

As in [LR], let $w_0$ be the image of $(I, -I)\in G\times G$ in $H$. It is given by
$$
w_1\pmatrix
I_n & 0\\
0 & -I_n\endpmatrix
w_1^{-1}=\pmatrix
0 & I_n\\
I_n & 0\endpmatrix,
\tag 5.6
$$
where $\ds w_1=\pmatrix I_n & I_n\\
I_n & -I_n\endpmatrix$, \qquad
$\ds w^{-1}_1=\pmatrix {1\over 2}I_n & {1\over 2}I_n\\  \\
{1\over 2}I_n & -{1\over 2}I_n\endpmatrix$, 
for the cases when $n$ is even. It is the long element of the Weyl group of $H$ modulo that of $M\simeq GL(V^d)$. 

We now define the intertwining operator 
$$
(M_{w_0}f)(x)=\int\limits_{U_P(F)}\ f(w_0u x)d u\quad (x\in H(F)),\tag5.7
$$
for every $f$ in $I(s,\chi)$ or to be exact, its space. To be more precise, one can set up a $K$--module isomorphism $f\longmapsto f_s : = f\cdot |\Delta|^s$ from $I(0,\chi)$ onto $I(s,\chi)$, where $K$ is a maximal compact subgroup of $H(F)$ such that $H(F)=P(F)K$. We may then write $M(s,\chi)f_s$ to denote $M_{w_0}f$. Note that
$$
M(s,\chi): I(s,\chi)\longrightarrow I(-s,\theta(\chi)^{-1}).\tag 5.8
$$
We recall that $M(s,\chi)$ converges absolutely for ${\text Re}(s)>>0$ and extends to a meromorphic function on all of $\bC$. Finally, assume $E=F$ and that $n/2$ is even if $G=O(n)$. Then the $w_0$--conjugate $\oP$ of $P$ will equal to opposite $P^-$ of $P$ which shares the same Levi as $P$. Define the map
$$j : I(-s,\chi^{-1})\longrightarrow\oI(s,\chi)\tag 5.9$$
by
$$jf(x)=f(w_0 x),\tag 5.10$$
where $\oI(s,\chi)$ is the representation induced from $\oP$. Then
$$j\cdot M(s,\chi)\  : \ I(s,\chi)\longrightarrow\oI(s,\chi)\tag5.11$$
equals the map $J_{\oP|P}$ of (4.15).

We now resume our generality and define the zeta--function that genalizes that of Godement--Jacquet [GJ].

Let $\pi$ be an irreducible admissible representation of $G(F)$ and $\tpi$ its contragredient. To avoid confusion in our notation, let ${p}={p}_\pi$ be the standard pairing
$${p}\ :\ \pi\otimes\tpi\longrightarrow \bC,\tag5.12$$
where $\pi$ and $\tpi$ are identified by their spaces. Given $\alpha\otimes\tal\in\pi\otimes\tpi$, consider the function ${p}(\pi(g)\alpha\otimes\tal)$, a matrix coefficient of $\pi$. We then set
$$Z(f_s,\alpha\otimes\tal)=\int\limits_{G(F)} f_s(i(g,1)){ p}(\pi(g)\alpha\otimes\tal)dg,\tag5.13$$
where $i : G(F)\times G(F)\hookrightarrow H(F)$ is the embedding as before. One can give the integration over $G^d(F)\backslash G(F)\times G(F)$ as in [LR] as it is more suitable for inductive arguments given in [LR] for unramified calculations as apposed to those in [GPSR]. 

The integral defining the zeta function $Z(f_s, \alpha\otimes\tal)$ converges for ${\text Re}(s)>> 0$ for a given $\pi$ and in fact for ${\text Re}(s)\ge 0$ if $\pi$ is square integrable (Lemma 2 of [LR]), and extends to a rational function in $q^{-s}$ on all of $\bC$. Moreover there exists a scalar--valued function $\Gamma(s,\pi,\theta)$ such that
$$Z(M(s)f_s, \alpha\otimes\tal)=\Gamma(s,\pi,\chi)Z(f_s,\alpha\otimes\tal).\tag 5.14$$
 This is Theorem 3 of [LR].
 
The point is that $\Gamma(s,\pi,\chi)$ is {\bf not} the $\gamma$--factor
$$\gamma(s,\pi\times\chi, \std,\psi)=\varepsilon(s,\pi\times\chi,\std,\psi){L(1-s,\tpi\times \chi^{-1}, \std)\over L(s,\pi\times\chi,\std)},\tag 5.15$$
where $\psi$ is a non-trivial (additive) character of $F$ and ''std'' denotes the standard representation of $\hat G$ which can be dropped from the notation.

To understand this discrepancy one has to compute the zeta function for the unramified data. As we will see the unramified zeta--function picks up an unwanted denominator, a phenomena  which also happens for the Rankin--Selberg $L$--functions for $GL_n\times GL_n$ and was remedied there by the introduction of an abelian Eisenstein series just as here [JS1]. (The same happens for exterior square $L$--functions for $GL_n$ in [JS2].)

We will follow [LR] again and refer to other papers as needed. We first give a table of the ''unwanted'' denominators. We will leave out the case $h=0$, which is needed as it will give the contribution of a global unitary group at places where it splits. Let
$$
d_H(s,\chi)=\cases
 L(s+\frac 12(n+1),\chi)\prod\limits^{n/2}_{j=1}\quad L(2s+2j-1,\chi^2)\\
\prod\limits^{n/2}_{j=1}\quad L(2s+2j-1, \chi^2)\\ \\
\prod\limits^{(n-1)/2}_{j=1}\quad L(2s+2j, \chi^2)\\ \\
\prod\limits^n_{j=1}\quad L\left(2s+j, \chi^0\,\eta^{n-j}_{E/F}\right),
\endcases\tag 5.16
$$
according as $h$ is symplectic, symmetric, $n$ even, symmetric, $n$ odd, or hermitian. In the last case $\chi^0 =\chi | F^*$ and $\eta_{_{E/F}}$ is the character of class field theory defining $E/F$. Finally all the $L$--functions are those of Hecke--Tate type for the field $F$. We now let $\alpha_0\in\pi$ and $\tal_0\in\tpi$ be $K$--fixed vectors with ${p}(\alpha_0\otimes\tal_0)=1$. Moreover, choose $f^0\in I(0,\chi)$ such that $f^0(k)=1$, for all $k\in K$.

The following is Proposition 3, pg. 333, of [LR], announced as Theorem 1.2 in\newline
[PSR].

\medskip\noindent
{\bf(5.17)}\ {\bf Proposition}.
{\it  Let $\alpha_0,\tal_0$ and $f^0_s$ be the $K$--fixed data just introduced. Moreover assume $\chi$ is unramified. Then
$$Z(f^0_s,\alpha_0\otimes\tal_0)=L(s+\frac 12, \pi\times\chi)\, /\, d_H(s,\chi).$$}

This proposition can be proved either following the arguments in pages 37--48 of [GPSR], (beware of typos!), which was generalized and corrected by Jian--Shu Li in [L], or the inductive arguments given in [LR] by Lapid and Rallis which were later generalized by Yamana to inner forms of unitary groups [Y1,2].

The work in [GPSR] then uses  a normalized Eisenstein series $E_H$, defined in pg. 32 of [GPSR], using the global version of $d_H$, to prove the finiteness of the poles of global standard $L$--functions.

\medskip\noindent
(5.18)\ {\bf Remark}.
The shift $s+\frac 12$ that appears in the $L$--function shows up in every case, including the inner forms, as well as the case of standard $L$--functions for $GL_n(D)$, if it is treated by means of doubling when $D$ is a central simple algebra over $F$. We will discuss this and its connection with basic functions later.

To obtain the $\gamma$--factor $\gamma(s,\pi\times \chi,\psi)$ attached to the standard $L$--function $L(s,\pi\times\chi)$ (see (5.15)) one needs to correct $\Gamma(s,\pi,\chi)$. This will eventually be given by a formal normalization of $M(s,\chi)$. But it is instructive to see what happens in the unramified setting. Without loss of generality we may assume $\chi=1$. To conform with the notation in [PSR], we introduce 
$$a_H(s)=a_H(s,1)=d_H(s-\langle \rho,\varepsilon_1\rangle, 1)\tag 5.19$$
when $n$ is even and we are in the symmetric case with $\langle \rho,\varepsilon_1\rangle=(n-1)/2$. We refer to [PSR] as well as equation (6.29) here for the symplectic case. Write
$$M(s)f^0_s=m(s)f^0_{-s},\tag 5.20$$
where $f^0_{-s}$ is the normalized spherical function in $I(-s, 1)$ and $m(s)$ is a scalar. Then the calculation in [LR] shows that
$$m(s)=a_H(s)\,/\,d_H(s),\tag 5.21$$
where $d_H(s):= d_H(s,1)$. Applying (5.14) to the unramified case one gets
$$m(s)Z(f^0_{-s},\alpha_0\otimes\tal_0)=\Gamma(s,\pi, 1) Z(f^0_s,\alpha\otimes\tal_0).\tag 5.22$$

Now Proposition (5.17) implies
$$m(s)d_H(-s)^{-1}L(-s+\frac 12,\tpi)=\Gamma(s,\pi, 1)L(s+\frac 12,\pi)d_H(s)^{-1}.\tag 5.23$$

In the unramified case
$$\aligned
\gamma(s+\frac 12,\pi,\psi): &= L(1-(s+\frac 12),\tpi) / L(s+\frac 12,\pi)\\
&=\Gamma(s,\pi,1){d_H(-s)\over d_H(s)}\,m(s)^{-1}\\ \\
&=\Gamma(s,\pi,1){d_H(-s)\over d_H(s)}{d_H(s)\over a_H(s)}\\ \\
&=\Gamma(s,\pi,1){d_H(-s)\over a_H(s)}.
\\
\endaligned\tag5.24
$$
We conclude that the correction factor
 $$\eta(s)=d_H(-s)/ a_H(s)\tag 5.25$$
will give the (correct) $\gamma$--factor as
$$\gamma(s,\pi,\psi)=\Gamma(s-\frac 12,\pi,1){d_H(\frac 12-s)\over a_H(s-\frac 12)}\tag 5.26$$
exactly as in Theorem 3.2 of [PSR] since it is identical to its definition in the line under equation [3--5] of the same theorem. This correction factor $\eta(s)$ is in fact a special case of a normalizing factor of $M(s,\chi)$ which we will explain after this; but let us record the unramified case as

\medskip\noindent
{\bf (5.27)}\ {\bf Proposition}
{\it The unramified $\gamma$--factor
$$\gamma(s,\pi,\psi):=L(1-s,\tpi)/ L(s,\pi)$$
is equal to
$$\gamma(s,\pi,\psi)=\Gamma(s-\frac 12,\pi,1)\eta(s-\frac 12),$$
where
$$\eta(s)=d_H(-s)/a_H(s)$$
is the corresponding correction factor.}
 
Note that 
$$m(s)\eta(s)d_H(s)\,/\, d_H(-s)=1\tag 5.28$$
which justifies the introduction of $M(s)$ and its normalizing factor $\eta(s)$. It simply removes the unwanted denominators $d_H(s)$ and $d_H(-s)$ from the zeta functions.

What generalizes this to ramified cases is an extension of definition of local coefficients of [Sh1] to our degenerate induced representation $I(s,\chi)$. To simplify the exposition, we will limit ourselves to the cases where $E=F$. We refer to [LR] for the general case. Then the $F$--points of the unipotent radical $U_P$ of $P$ is isomorphic to
$$X_n(\varepsilon)=\{X\in M_n(F)|\, X +\varepsilon^t X=0\},\tag5.29$$
where $\varepsilon=1$ or $-1$, according to whether $h$ is symmetric or symplectic, respectively, through the exponential map. We now fix $A\in X_n(\varepsilon)$ with $\det(A)\ne 0$. We can then define a character of $U_P(F)$ by
$$\psi_A\,:\, T\longmapsto \psi(tr(TA))\tag 5.30$$
for all $T\in X_n(\varepsilon)$. 

Now consider the induced representation $I(\psi_A)=\Ind^{H(F)}_{U_P(F)}\psi_A$. Then one knows that
$$\dim(\Hom_{H(F)}(I(s,\chi), I(\psi_A)))\le 1\tag 5.31$$
and indeed equal to 1 since $\det(A)\ne 0$ (cf. [Ka]). In fact, the functional
$$\ell_{\psi_A}(f)=\int\limits_{U_P(F)} f(w_0 u)\overline{\psi_A}(u)\, du\tag 5.32$$ 
is a non-zero candidate for this space. Here $f\in I(s,\chi)$.

One can define a similar functional $\ell'_{\psi_A}$ for $I(-s,\chi^{-1})$. A {\bf degenerate local coefficient} $c(s,\chi, A, \psi)$ is now defined by means of multiplicity one in (5.31) through
$$\ell'_{\psi_A}(M_{w_0}(f))=c(s,\chi,A,\psi)\ell_{\psi_A}(f).\tag 5.33$$
We now normalize $M_{w_0}(s,\chi)=M_{w_0}$ by $c(s,\chi,A,\psi)$ to get:
$$M_{w_0}^\ast(s,\chi,A,\psi)=c(s,\chi,A,\psi)^{-1}M_{w_0}(s,\chi).\tag 5.34$$

We note that as in equation (16), page 327, of [LR]
$$M_{w_0}^\ast(-s,\chi^{-1},A,\psi)M_{w_0}^\ast(s,\chi,A,\psi)=I,\tag 5.35$$
i.e., $M_{w_0}^\ast(s,\chi,A,\psi)$ is normalized.

With $\Gamma(s,\pi,\chi)$ as in (5.14) one now defines
$$\Gamma(s,\pi,\chi,A,\psi)=\Gamma(s,\pi,\chi)\ c(s,\chi,A,\psi)^{-1}.\tag 5.36$$
In the unramified case and after choosing $\psi$ and $A$ suitably, one gets
$$\theta(s) c(s,1,A,\psi)^{-1}=\eta(s)\tag 5.37$$
with $\eta(s)$ as in (5.25), where $\theta(s)=1$ if $\eta$ is symmetric, while $\theta(s)=\gamma(s+\frac 12,\psi)$, the Hecke--Tate $\gamma$--factor, if $h$ is symplectic.

The $\gamma$--factors $\gamma(s,\pi\times\chi,\psi)$ is then defined by equation (25), page 337, of [LR] for any irreducible admissible representation $\pi$ and character $\chi$ of $E^*$. It does not depend on the choice of $A$.

\medskip\noindent
(5.38)\ {\bf Remark}. 
As in the case of local coefficients, there is a global analogue of $\ell_{_{\psi_A}}$ when one takes the $\psi_A$--Fourier coefficient of the corresponding (degenerate) Eisenstein series which unfolds to a product of local functional on local induced representations. Using the functional equation of the Eisenstein series, equation (5.33) and finally local functional equation (5.14), this leads to the  global functional equation for $L(s,\pi\times\chi)$ (cf. equation (24) and discussion in page 340 of [LR]). The steps are clearly parallel to those of local coefficients and its global theory in the generic case, or more generally when the representation has other models [FG].

\bigskip\noindent
{\bf 6.\ The basic function}

We shall now determine the basic function $\bL^{\std}(s)(g)$ attached to the standard $L$--function $L(s,\pi)$ for classical groups (cf. Section 3), and show that it is fixed by the Fourier transform $\{\sF_{Q|P}\}$ on $\sS(X_P)$, the space of Schwartz functions defined in Section 4.

This can be done at the level of doubling method, realizing the Schwartz space through the spaces of induced representations by means of their Mellin transform (4.6), or through the Schwartz space of Braverman--Kazhdan itself [BK2].

We recall from the discussions in Sections 3 and 4 of this paper, as well as theorem 7.4.9 and other discussions in Section 7 of [Li], that the reductive monoid $X\ :\ =X_{\std}$ attached to the standard representation of $\hat G$ contains $X_P$ and in fact $X=X_P\sqcup\{0\}$, the affine closure of $X_P$. We thus have
$$X^+\underset{\text open}\to\subset\, X_P \underset{\text open}\to\subset\,X\tag 6.1$$
with $X^+$ the main orbit of $X_P$ under the action of $M_{ab}\times G\times G$ on it (cf. [BK2,Li]); it is open and dense in $X_P$ and isomorphic to $M_{ab}\times G$. Moreover, $X_P$ is quasi--affine whose affine closure is $X$.

With $M_{ab}\times G$ as the (reductive)  group of units of $X$, the function $c$ of Section 3 is now
$$\gathered
M(F)    \buildrel{\Delta}\over\longrightarrow M_{ab}(F)\\
m\longmapsto \overline m.
\endgathered\tag 6.2
$$
where $\Delta=\det$ under the isomorphism $M\simeq GL(V^d)$. More precisely,
$$c(\overline m, g)=\Delta(m),\tag 6.3$$
$\overline m\in M_{ab}(F)$ and $g\in G(F)$.

To connect our zeta functions to Braverman--Kazhdan's generalization of Godement--Jacquet, we again resort to the discussions in Section 4. Using the reductive group of units of $X$, $M_{ab}\times G$, we can now start with a function $\xi\in C^\infty_c(X_P(F))$ and define a function in $I(s,\chi)$ by
$$x\longmapsto\int\limits_{M_{ab}(F)}\chi_s(m)\xi((\overline m)^{-1}\cdot x)\, d\overline m,\tag 6.4$$
with $\delta_P(m)$ built into the measure (cf. the discussion just before Remark 7.1.2 in [Li]), where $m$ is the preimage of $\overline m$ under
$$P / P_{\der}\simeq M / M_{\der}\simeq M_{ab}.$$

Let $\sM_\pi(G(F))$ denote the space of matrix coefficients of $\pi$. Integrating functions in $C^\infty_c(X_P(F))\otimes\sM_\pi (G(F))$ over $M_{ab}(F)\times G(F)$ will now give our zeta functions.

We will first give a treatment within the doubling method. We recall from Proposition 5.17 that
$$Z(f^0_s,\alpha_0\otimes\tal_0)=L(s+\frac 12,\pi)d_H(s)^{-1}.\tag 6.5$$
To get the basic function we need to correct the shift $\frac 12$. We therefore define our basic function as
$$\bL^{\std}(s)(g)=d_H(s-\frac 12)f^0_{s-\frac 12}(g).\tag 6.6$$

We now calculate the Fourier transform in the doubling setting of $f^0_{s-\frac 12}(g)d_H(s-\frac 12)$, i.e., 
$$M^*(f^0_{s-\frac 12}(g))d_H(s-\frac 12).\tag 6.7$$
Using (5.28) we conclude that (6.7) equals to
$$\aligned
\eta(s-\frac 12)m(s-\frac 12)&f^0_{\frac 12 -s}(g) d_H(s-\frac 12)\\
&=f^0_{\frac 12 -s}(g) d_H(\frac 12 -s).\endaligned\tag 6.8$$
To agree with the treatment in [BK2], we need further to go to the setting of opposite parabolic $\oP$. The map $j$, defined by (5.10), can now be applied as 
$$j\,:\, I(\frac 12 -s,{\bold 1})\longrightarrow \oI(s-\frac 12,{\bold 1}), \tag 6.9$$
where $\oI(t,{\bold 1}), t\in \bC$, is the representation of $H(F)$ induced from $\oP=M\oN$, i.e.,
$$\oI(t,{\bold 1})=\Ind^{H(F)}_{M(F)\oN(F)}\, \bold{1}_t\otimes{\bold 1}.\tag 6.10$$ 

The equation (6.8) now changes to
$${\bar f}^0_{s-\frac 12}\, (g) \overline{d_H}(s-\frac 12)\tag 6.11$$
and is thus preserved by the Fourier transform $j.M^*$ applied to it. Here $\overline f^0_{s-\frac 12}$ is the normalized unramified function in $\oI(s-\frac 12,{\bold 1})$. We note that $\overline d_H(t)$, defined using $\oI(t,\bold 1)$, where $\oP$ is used for the doubling method instead of $P$, $t\in \bC$, satisfies 
$$\overline{d_H}(t)=d_H(-t).\tag 6.12$$
We record this as

\medskip\noindent
{\bf (6.13)}\ {\bf Proposition}
{\it The basic function 
$$\bL^{\std} (s)(\overline m,g)=d_H(s-\frac 12)f^0_{s-\frac 12} (g),$$
whose integral against ${p}(\pi(g)\alpha_0\otimes\tal_0)$ gives $L(s,\pi)$, is fixed by the Fourier transform $j.M^*$.} 

We shall now address the problem at the level of Schwartz functions (cf. Section 4). The operator $\sF_{P|\sQ}$, which will be shortly defined, identifies $\sS(X_Q)\subset L^2(X_{\sQ})$ with $\sS(X_P)\subset L^2(X_P)$. We will use $\sS(M)=\sS(H, M)$ to denote these isomorphic spaces and thus
$$\sS(M):=\sS(X_P)=\sum\limits_{\sQ}\sF_{P|\sQ}(C^\infty_c(X_{\sQ})).\tag6.14$$

To have complete agreement with generalized Godement--Jacquet theory, we need a function on $M_{ab}(F)\times G(F)$, the $F$--points of the group $X^+$ of units of our reductive monoid $X=X_{\std}$, which when integrated against ${p}(\pi (g)\alpha_0\otimes\tal_0)$ over $M_{ab}(F)\times G(F)=X^+(F)$ gives $L(s,\pi)=L(s,\pi,\std)$.
In our situation the only possibilities are $P$ and $\oP$, the opposite parabolic.

Choose $\xi^0_P\in \sS(M)^K$, $K=H(O_F)$, considered as a smooth function on $X^+(F)$ and $X^+(O_F)$--invariant such that 
$$\int\limits_{M_{ab}(F)}|c(\om,g)|^{s-\frac 12}\xi^0_P((\om)^{-1} ,g)d{\om}= d^P_H(s-\frac 12)f^0_{s-\frac 12}(i(g,1)).\tag 6.15$$
In view of Conjecture 7.1.5 of [Li], we may assume the integral in (6.15) is convergent for $s$ in an appropriate cone. (This is not automatic since $\xi^0_P$ is not of compact support in $H(F)$ and its restriction to $X^+(F)$ may not be as such neither.) We will show later that such $\xi^0_P$ exists.

Define the {\bf basic function}
$$\bL^{\std}_P(s)(\om,g)=|c(\om,g)|^{s-\frac 12}\xi^0_P((\om)^{-1},g).\tag 6.16$$
Then
$$\int\limits_{M_{ab}(F)\times G(F)} \bL^{\std}_P (s)(\om,g) p(\pi(g)\alpha_0\otimes\tal_0)d\om dg=L(s,\pi,\std),\tag 6.17$$
justifying the name.

To continue we need to appeal to the discussion of the Braverman--Kazhdan's paper [BK2] which we discussed partly in Section 4.

We recall the representation $L$ of $\hT=Z(\hM)$, introduced in Section 4, which was used to define $\sF_{Q|P,\psi}$ and specialize it to the cases of interest for us. Let $P$ and $Q$ be parabolic subgroups of $H$, sharing the same Levi subgroup $M$. Thus $P=M U_P$ and $Q=M U_Q$. Let $\hat P, \hat Q, \hat M, {\hat {U}}_P$ and ${\hat{U}}_Q$ be the dual groups. Let ${\hat{\frak p}}=Lie(\hat{P})$, ${\hat{\frak q}}=Lie(\hat Q)$, ${\hat{\frak u}}_{\frak p}=Lie({\hat{U}}_P)$ and ${\hat{\frak u}}_{\frak q}=Lie({\hat{U}}_Q)$.
Set $\hat{\frak u}_{\frak p,\frak q}={\hat{\frak u}}_{\frak p}/{\hat{\frak u}}_{\frak p}\cap {\hat{\frak u}}_{\frak q}$. 
Let $\{e,h,f\}$ be a principal (regular) $SL_2(\bC)$--triple in $\hat{\frak m}=Lie(\hM)$. The adjoint action of $\hM$ on ${\hat{\frak u}}_{\frak p,\frak q}$ restricts to a representation of this $SL_2(\bC)$--triple.
Let $({\hat{\frak u}}_{\frak p,\frak q})^e$ be the set of highest weight vectors for $e$ in ${\hat{\frak u}}_{\frak p,\frak q}$. With notation as in Section 4 (cf. [BK2]), we let $L=({\hat{\frak u}}_{\frak p,\frak q})^e$. In our setting $\hat{\frak u}_{\frak p}\cap\hat{\frak u}_{\frak q}=\{0\}$ and thus ${\hat{\frak u}}_{\frak p,\frak q}={\hat{\frak u}}_{\frak p}$. 

We will now complete the definition of the Fourier transform of Braverman--Kazhdan [BK2] which we addressed in Section 4. With $L=({\hat{\frak u}}_{\frak p,\frak q})^e$, we define the distribution $\eta_{_{L,\psi}}=\eta_{_{P,Q,\psi}}$ as in (4.23) and set
$$\sF_{P|Q,\psi}\, :\, \sS(X_Q)\longrightarrow\sS(X_P)$$
by
$$\sF_{P|Q,\psi}=\eta_{_{P,Q,\psi}}\cdot\sR_{P|Q},$$
where $\sR_{P|Q}$ is as in equation (4.14).

Let
$$\kappa: SL_2(\bC)\longrightarrow\hM\tag 6.18$$
be the homomorphism attached to our triple. For $t\in\bC^*$, let $H_t=\kappa(\pmatrix t & 0\\
0 & t^{-1}\endpmatrix)$. Then by Jacobson--Morozov, [C], page 139, $H_{q^{\frac 12}}$ gives the Satake parameter for the trivial representation of $M(F)$. It follows that the adjoint action of $H_{q^{\frac 12}}$ will be given by multiplication by $H_q$ in a basis of root vectors in $\hat{\frak u}_{\frak p}$. 

Changing $s$ to $s+{n\mp 1\over 2}$ according as $G=O(n)$ or $Sp(n)$, so as to get the normalized induction, one will have
$$\chi_s=\sum\limits^n_{i=1}\{s-[\frac 12(n+1)-i]\} x_i\tag 6.19$$
in both cases, for now ``normalized'' inducing data of [GPSR] (cf. Lemma 5.2 of [GPSR]). We recall the normalized induction from the characters
$$\mu_i(t_i)=|t_i|^{\frac 12(n+1)-i}$$
giving the trivial representation of $M(F)=GL_n(F)$ and its Satake parameter.

The character $\chi_s$ of $M(F)$ may be regarded as one of $M_{ab}(F)$. Let $T_H$ be the standard maximal torus of $H$ given by coordinate functions $x_i$ and contained in $M$. Then using the map
$$T_H\subset M\longrightarrow M/M_{\der}:=M_{ab},\tag 6.20$$
one can lift $\chi_s$ to a character of $T_H(F)$. The character $\chi_s$ in (6.19) may be considered to be this lift, and in what follows sometimes denoted by $\widetilde\chi_s$.

With notation as in Theorem 5.10 of [BK2], define a function on $Z(\hM)$ by
$$d_P(z)=\det(1-q^{-1}H^{-1}_q\cdot z)|(\hat{\frak u}_{\frak p})^e\qquad\quad (z\in Z(\hM)).\tag 6.21$$
The variable $z$ acts by adjoint action. This function and its ``dual'' will define the Fourier transform on the canonical basis of $K$--invariant functions in $\sS(M)$ as explained in Lemma 5.11 of [BK2]. We shall now try to relate them to our function $d_H$, normalizing factors and normalized operators.

We resort to the calculations in [GPSR] which will now be normalized by sending $s$ in [GPSR] to $s+{n\mp 1\over 2}$ according as $G=O(n)$ or $Sp(n)$, respectively. When $G=O(n)$, the intertwining scalar $m(s)$, defined in equations (5.20) and (5.21) can now be written as
$$\aligned
m(s)&=\prod\limits^n_{\ell=1}\zeta((\chi_s, x_\ell+x_{\ell+1}))/\zeta(1+(\chi_s, x_\ell+x_n)) \\
\\
&=\prod\limits^{n/2}_{\ell=1}\zeta((\chi_s,x_\ell+x_{\ell+1}))\bigg/ \prod\limits^{n-1}\Sb{\ell=1}\\ \ell\equiv 1(2)\endSb \zeta(1+(\chi_s, x_\ell+x_n))
\\
\endaligned\tag6.22
$$
with the standard $\zeta$--function for the field $F$. Here $x_i$, $1\le i\le n$, $n$ even, are the coordinate characters of the maximal torus of $SO(2n,\bC)\subset O(2n,\bC)$ and contained in $\hM=GL_n(\bC)$, giving the roots in $\hat{\frak u}_{\frak p}$.

This is the function $c_{w_0}(s)$ in page 29 of [GPSR]. Moreover, per formulas in [GPSR], we are using $(\chi_s, x_\ell+x_m)$ to denote $(\chi_s,\alpha^\vee)$, where $\alpha=x_\ell+x_m$ since $(x_\ell+x_m)^\vee=x_\ell+x_m$, $\ell\ne m$.

The case of $G=Sp(n)$ is similar. One needs to note that $(2x_i)^\vee=x_i$ which gives the extra quotient of $\zeta$--function coming from $(2x_i)^\vee=x_i$ as follows
$$
\prod\limits^n_{i=1}{\zeta((\chi_s, x_i))\over\zeta(1+(\chi_s, x_i))}=\prod\limits^n_{i=1}{\zeta(s+{n+1\over 2}+i-n-1)\over\zeta(s+{n+1\over 2}+i-n)}\tag6.23$$
which simplifies to
$$\zeta(s-{n-1\over 2})/\zeta(s+{n+1\over 2})
\tag6.24$$
as a factor in $c_{w_0}(s)$. The denominator in (6.24) is precisely
$$\aligned
\zeta(s+\frac 12(n+1))&=\zeta(1+(\chi_s, x_n))\\
&=\beta(s)\endaligned\tag 6.25$$
of page 4590 of [PSR] as a factor of $d_H(s)$, while the numerator of (6.24) is exactly
$$\aligned
\zeta(s-\frac 12(n-1))&=\zeta((\chi_s,x_1))\\
&=\tbe(s)\endaligned\tag 6.26$$
introduced in page 4591 of [PSR] as a factor of $a_H(s)$.

Thus for $G=Sp(n)$ we get the formula
$$m(s)=\prod\limits^{n-1}_{\ell=1}{\zeta((\chi_s,x_\ell+x_{\ell+1}))\over\zeta(1+(\chi_s, x_\ell+x_n))}\cdot{\zeta((\chi_s,x_1))\over \zeta(1+(\chi_s, x_n))}\tag 6.27$$
using formula (5.21).

We now express the normalizing factor $\eta(s)$, introduced by equation (5.25), as a product of $\gamma$--functions as suggested in Section 4. In the case of double covering of $Sp(n)$ this is proved by Gan in [Gan]. Recall the normalizing factor
$$\eta(s)=d_H(-s)/ a_H(s).\tag 6.28$$
In the new setting the function $a_H$(s), defined in [PSR] (also see equation (5.19)), is given by
$$a_H(s)=\prod\limits^{n/2}_{\ell=1}\zeta((\chi_s,x_\ell+x_{\ell+1}))\cdot
\cases 1 & G=O(n)\\ \zeta((\chi_s,x_1)) & G=Sp(n).\endcases\tag 6.29$$
Moreover
$$d_H(-s)=\prod\limits^{n/2}_{\ell=1}\zeta(1-(\chi_s,x_\ell+x_{\ell+1}))\cdot
\cases 1 & G=O(n)\\ \zeta(1-(\chi_s,x_1)) & G=Sp(n),\endcases\tag 6.30$$
since $\zeta((\chi_{-s}, x_n)+1)=\zeta(1-(\chi_s, x_1))$ for $G=Sp(n)$.

Thus
$$
\eta(s)=\prod\limits^{n/2}_{\ell=1}{\zeta(1-(\chi_s,x_\ell+x_{\ell+1}))\over\zeta((\chi_s, x_\ell+x_{\ell+1}))}\cdot{\zeta(1-(\chi_s,x_1))\over \zeta((\chi_s, x_1))}.\tag6.31$$
when $G=Sp(n)$, while
$$\eta(s)=\prod\limits^{n/2}_{\ell=1}{\zeta(1-(\chi_s,x_\ell+x_{\ell+1}))\over\zeta((\chi_s, x_\ell+x_{\ell+1}))}\tag 6.32$$
if $G=O(n)$.

Given a variable $T$, define
$$\gamma(T)=\zeta(1-T)/\zeta(T).\tag 6.33$$
We have proved:

\medskip\noindent
{\bf (6.34)}\ {\bf Proposition}
{\it The normalizing factor $\eta(s)$ is a product of $\gamma$--functions. More precisely}
$$\eta(s)=\prod\limits^{n/2}_{\ell=1}\gamma((\chi_s,x_\ell+x_{\ell+1}))\cdot
\cases \gamma((\chi_s,x_1)) & G=Sp(n)\\ 1 & G=O(n).\endcases$$

Finally, we reformulate the normalized operators as follows:
$$\aligned
m(s)\eta(s)&={{\prod\limits^{n/2}_{\ell=1}\zeta((\chi_s,x_\ell+x_{\ell+1}))\over\prod\limits^{n-1}\Sb{\ell=1}\\ \ell\equiv 1(2)\endSb \zeta(1+(\chi_s, x_\ell+x_n))}\cdot{\zeta((\chi_s, x_1))\over\zeta(1+(\chi_s, x_n))}}\\
\\
&\cdot \prod\limits^{n/2}_{\ell=1}{\zeta(1-(\chi_s,x_\ell+x_{\ell+1}))\over \zeta((\chi_s, x_\ell+x_{\ell+1}))}\cdot{\zeta(1-(\chi_s, x_1))\over\zeta((\chi_s, x_1))}\\
\\
&={\prod\limits^{n/2}_{\ell=1}\zeta(1-(\chi_s,x_\ell+x_{\ell+1}))\, \zeta(1-(\chi_s,x_1)) \over\prod\limits^{n-1}\Sb{\ell=1}\\ \ell\equiv 1(2)\endSb \zeta(1+(\chi_s, x_\ell+x_n))\, \zeta(1+(\chi_s, x_n))}\endaligned\tag 6.35$$
when $G=Sp(n)$, with the factor
$$\zeta(1-(\chi_s,x_1))/\zeta(1+(\chi_s,x_n))\tag6.36$$
missing when $G=O(n)$.

The following lemma determines the set $({\hat{\frak u}}_{\frak p})^e$ of highest weights of $e$ of our triple $(e,h,f)$ in ${\hat{\frak u}}_{\frak p}$.

\medskip\noindent
(6.37)\ {\bf Lemma}.
{\it The set $({\hat{\frak u}}_{\frak p})^e$ consists of root vectors in ${\hat{\frak u}}_{\frak p}$ attached to roots $x_\ell+x_{\ell+1}$, $1\le\ell\le\frac n2$ for $G=O(n)$, and the same set together with $x_1$ if $G=Sp(n)$.}

\medskip\noindent
{\bf Proof}. 
Element $e$ of the triple $\{e,h,f\}$, being regular unipotent, can be represented by a sum of root vectors $X_\alpha$, $\alpha=x_i-x_{i+1}$, $1\le i\le n-1$, simple roots of $GL_n(\bC)$. Assume $G=O(n)$, $n$ even. The root spaces in ${\hat{\frak u}}_{\frak p}$ add up to a direct sum on which $e$ acts by adjoint action. It will act irreducibly on
$${\hat{\frak u}}_{{\frak p},1}=\bigoplus\limits^n_{\ell=2}\bC X_{x_1+x_\ell},$$
where $\bC X_\alpha$ gives the root space of given root $\alpha$. We observe that ${\hat{\frak u}}_{{\frak p},1}$ has $X_{x_1+x_2}$ and $X_{x_1+x_n}$ as the highest and lowest weight vectors, respectively.

We now consider ${\hat{\frak u}}_{\frak p}/{\hat{\frak u}}_{{\frak p},1}$ and the image ${\hat{\frak u}}_{{\frak p},2}$ of $\bigoplus\limits^n_{\ell=3}\bC X_{x_2+x_\ell}$ in it. The element $e$ again acts irreducibly with $X_{x_2+x_3}$ and $X_{x_2+x_n}$ as the highest and lowest weight vectors, respectively. We continue in this way up to ${\hat{\frak u}}_{{\frak p},n/2}$.

To proceed, we note that the $SL_2(\bC)$--triple $\{e,h,f\}$ defines a non-trivial Weyl group representative $w$ by the standard formula $w=exp (e)exp(-f)exp(e)$ which may be considered as a representative for the long element of the Weyl group of $\hM=GL_n(\bC)$ by regularity of the triple. It is clear that $\kappa(w)$ fixes the highest weight attached to the root $x_{n/2}+x_{(n/2)+1}$, while sending the one attached to $x_{n-\ell}+x_{n-\ell+1}$, to that of root $x_\ell+x_{\ell+1}$, $1\le\ell< n/2$. Since weights of irreducible representations of $SL_2(\bC)$ are Weyl group conjugate, the ones for $x_\ell+x_{\ell+1}$, $\ell> n/2$, cannot be highest weights for $e$.

For $G=Sp(n)$, $x_1$ and $x_n$ become the roots giving the highest and lowest weights of corresponding extra representation that appear in this case.

The map (6.20) also leads to a map $\Lambda_*(T_H)\longrightarrow \Lambda_*(M_{ab})$. The cocharacter $\lambda_i\in\Lambda_*(M_{ab})$ of Section 4 for the action of $\hM_{ab}$ on $L=(\hat{\frak u}_{\frak p})^e$ can be lifted to cocharacter $\widetilde \lambda_i\in\Lambda_*(T_H)$. Moreover for every character $\theta$ of $M_{ab}(F)$, let $\widetilde\theta$ be a lift to $T_H(F)$ via map (6.20). Then $\theta\cdot\lambda_i=\widetilde\theta\cdot\widetilde\lambda_i$.

To calculate $\gamma$--functions (4.20) in Proposition 4.27, we need to determine half--eigenvalues of $h$ on $L=(\hat{\frak u}_{\frak p})^e$, where $h$ is the semisimple element in our regular $sl_2(\bC)$--triple $\{e,h,g\}$. By Jacobson--Morozov's theorem [C], we can write
$$(\frac 12 h,\tlam_i)=(\sum\limits^n_{j=1} s_j x_j, \tlam_i),$$
$\tlam_i=x_i+x_{i+1}$, $1\le i\le n/2$, where $s_j={n+1\over 2}-j$, $1\le j\le n$, with the choices of $\tlam_i$ given by Lemma 6.37. 
Let $\tchi'_s$, $s\in\bC$, be the character $|\det(\cdot)|^s=|\cdot|^{s\sum\limits^n_{j=1} x_j}$. The $\gamma$--functions (4.20) in Proposition 4.27, $\gamma(s_i,\tchi'_s\cdot\tlam_i)$, are now
$$\gamma(s_i,\tchi'_s\cdot\tlam_i)=\int\limits_{F^*} |t|^{(\sum\limits^n_{j=1}(s+ s_j)x_j,\tlam_i)}\psi(t)dt.$$
This is evidently equal to $\gamma((\tchi_s,x_{n-i}+x_{n-i+1}))$, where $\tchi_s$ is given by (6.19). We therefore have the following

\medskip\noindent
(6.38)\ {\bf Corollary}.
{\it The normalizing factor $\eta(s)$, given explicitly by equations (6.31) and  (6.32), is the same as the one defined by $M(\eta_{_{L,\psi}})$ of Proposition (4.27) by Braverman--Kazhdan for the data $({\hat{\frak u}}_{\frak p})^e$, upon realizing $\tlam_i=x_i+x_{i+1}\in\Lambda_*(T_H)$, $1\le i\le n/2$, where $G=O(n)$, together with $\tlam_0 :=x_1$, if $G=Sp(n)$.}

Finally, we like to discuss the normalized operator $\sF_{P|Q}:=\sF_{P|Q, \psi}$ of Braverman--Kazhdan [BK2], their Fourier transform, and its effect on our basic function (6.16). We recall the function $d_P(z)$ of [BK2], $z\in Z(\hM)$, defined by equation (6.21) here. We have

\medskip\noindent
(6.39)\ {\bf Lemma}.
{\it a)\ The numerator of the normalized operator $m(s)\eta(s)$ equals $d_P(q^{s})^{-1}$. More precisely, 
$$\aligned
d_P(q^{s})^{-1}&=\prod\limits^{n/2}_{\ell=1}\zeta(1-(\chi_s,x_\ell+x_{\ell+1}))\cdot
\cases 1 & G=O(n)\\ \zeta(1-(\chi_s,x_1)) & G=Sp(n)\endcases\\
&=d_H(-s).\endaligned$$
b)\ The denominator of $m(s)\eta(s)$ equals $d_{\oP}(q^{-s})^{-1}$. More precisely, 
$$\aligned
d_{\oP}(q^{-s})^{-1}&=\prod\limits^{n-1}\Sb{\ell=1}\\ \ell\equiv 1(2)\endSb\zeta(1+(\chi_s,x_\ell+x_n))\cdot
\cases 1 & G=O(n)\\ \zeta(1+(\chi_s,x_n)) & G=Sp(n)\endcases\\
&=d_H(s).\endaligned$$
Here we realize $Z(\hM)\simeq\bC^\ast$, $\hM=GL_n(\bC)$.}

\medskip\noindent
{\bf Proof}. Part a) can be proved by direct calculations. It should also follow the fact that $H_q$ represents the adjoint action of the Satake parameter of the trivial representation of $GL_n(\bC)$. Part b) should be a consequence of the fact that $d_{\oP}$ is the dual of $d_P$ and thus given by lowest weights of the action of $e$. Details are left to the reader.

We will now show that the Fourier transform $\sF_{\oP|P}$ fixes our basic function (6.16) defined on $X^+(F)$. Let $\xi^0_P\in\sS(M)^K$ be as in (6.15), i.e., the function whose restriction to $X^+(F)$ was used to define the basic function $\bL^{\std}_P(s)$ which we shall now explicate. As explained in Section 5 of [BK2], as well as Section 3 of [BK3], the cocharacter lattice $\Lambda_*(M_{ab})$ of $M_{ab}$ parameterizes a basis for $\sS(M)^K$, $K=H(O)$. More precisely, as we discussed earlier, the map
$$T_H\subset M\longrightarrow M/M_{\der}:=M_{ab}\tag 6.40$$
leads to the restriction map $\Lambda_*(T_H)\longrightarrow \Lambda_*(M_{ab})$, where $T_H$ is the maximal torus of $H$ contained in $M$ fixed earlier. Given $\gamma\in\Lambda_*(M_{ab})$, let $\tga$ be any lift of $\gamma$ to $\Lambda_*(T_H)$. Consider $X_P^{\tga}$, the $K$--orbit of $\tga({\varpi}_F)\mod P_{\der}$. It only depends on $\gamma$ and thus we set $X^\gamma_P : = X^{\tga}_P$. For each $\gamma\in\Lambda_*(M_{ab})$, one defines 
$$\delta_{P,\gamma}(x)=\cases q^{\langle \gamma,\rho_P\rangle} &x\in X^\gamma_P\\
0 &{\text{otherwise}}.\endcases\tag 6.41
$$
In fact, $X_P=\cup_\gamma X^\gamma_P$, $\gamma\in\Lambda_*(M_{ab})$, and the functions $\delta_{P,\gamma}$ make a basis for $C^\infty_c(X_P(F))^K$. Then, as explained in Section 3.12 of [BK3] and Section 5 of [BK2], $\Lambda_*(M_{ab})$ acts on functions $\delta_{P,\gamma}$ by
$$\mu(\delta_{P,\gamma})=q^{\langle\mu,\rho_{_P}\rangle}\delta_{P,\gamma+\mu}.\tag 6.42$$
The inverses of polynomials $d_P(z)$ and $d_{\oP}(z)$ can be presented as elements in the symmetric algebra Sym$((\hat{\frak u}_{\frak p})^e)$ of $(\hat{\frak u}_{\frak p})^e$. Since $\bC[\Lambda_*(M_{ab})]\simeq\bC[\hM_{ab}]$, where $\bC[\Lambda_*(M_{ab})]$ is the group algebra of $\Lambda_*(M_{ab})$ and $\bC[\hM_{ab}]$ is the algebra of regular functions on $\hM_{ab}$, (cf. [BK2]), the coefficients of $d^{-1}_P$ and $d^{-1}_{\oP}$ are symmetric polynomials on $\Lambda_*(M_{ab})$.

The function $\xi^0_P\in\sS(M)^K$ in equation (6.15) can be taken
$$\xi^0_P(x)=d^{-1}_{\oP}(\delta_{P,0})(x).$$
$x\in X_P(F)$, where the action is according to (6.42) and $d^{-1}_{\oP}:=d^{-1}_{\oP}(1)$. We point out that
$$d^{-1}_{\oP}(\delta_{P,0})| X^\gamma_P=\Phi^L_{P,0}|X^\gamma_P,$$
where $\Phi^L_{P,\mu}=c_{_{P,\mu}}$, $\mu\in \Lambda_*(M_{ab})$, with $\Phi^L_{P,\mu}$ and $c_{_{P,\mu}}$ defined as in [BK2], (cf. equation (5.5) in [BK2]), with $L=\Sym((\hat{\frak u}_{\frak p})^e)$. We have

\medskip\noindent
(6.43)\ {\bf Lemma}.
{\it The Mellin transform of $\xi^0_P=\delta^{-1}_{\oP}(\delta_{P,0})$ equals 
$$\int\limits_{M_{ab}(F)} |c(m,g)|^s d^{-1}_{\oP}(\delta_{P,0})(m^{-1}g)dm=d_H(s)f^0_s(i(g,1)).$$}

\medskip\noindent
{\bf Proof}. 
The operator $d_{\oP}$ is a linear combination of $\tmu$ in the coroot lattice $\Lambda_*(T_H)$ which are weights of the adjoint action of ${\hT}_H$ in $({\hat{\frak u}}_{\frak p})^e$. Then $d^{-1}_{\oP}$ will be an infinite series on symmetric polynomials on these weights. Each $\tmu$ will act on $\delta_{P,0}$ by
$$\tmu(\delta_{P,0})(m^{-1}g)=\delta_{P,\mu}(m^{-1}g).$$
Let $\tmu$ be appearing in the expansion of $d^{-1}_{\oP}$. Then its contribution to the Mellin transform is
$$\int\limits_{M_{ab}(F)} |c(m,g)|^s\delta_{P,\tmu}(m^{-1}g) d m.$$
Changing $m$ to $m\tmu(\varpi_{_F})^{-1}$,  one gets
$$\int\limits_{M_{ab}(F)} |c(m\tmu(\varpi_{_F})^{-1},g)|^s\ \delta_{P,\tmu}(\tmu(\varpi_{_F})m^{-1}g)\  \delta_P(\tmu(\varpi_{_F}))^{-1} dm,$$
where $\delta_P$ is the modulus character of $P$. 
By definition
$$\delta_P(\tmu(\varpi_{_F}))^{-1\!/\!2}\delta_{P,\tmu}(\tmu(\varpi_{_F})m^{-1}g)=\delta_{P,0}(m^{-1}g).$$
Thus the contribution is 
$$\aligned
&\delta_P(\tmu(\varpi_{_F}))^{-1\!/\!2} |c(\tmu(\varpi_{_F}))|^{-s}\int\limits_{M_{ab}(F)}|c(m,g)|^s\delta_{P,0}(m^{-1}g)dm\\
&=\delta_P(\tmu(\varpi_{_F}))^{-1\!/\!2}|\mu(\varpi_{_F})|^{-s} f^0_s(g),
\endaligned
$$
where $f^0_s(g) : = f^0_s(i(g,1))$, with $i$ the embedding of $G\times G$ into $H$, which agrees with the embedding of $X^+=M_{ab}\times G\hookrightarrow X_P$. We remark that each term $\delta_P(\tmu(\varpi_{_F}))^{-1\!/\!2} | c(\tmu(\varpi_{_F}))|^{-s}$ accounts for the way elements in $X_*(M_{ab})$ act on $\delta_{P,\gamma}$ by equation (6.42). Now summing up over all the contribution  we get that our Mellin transform equals 
$$d_{\oP}(q^{-s})^{-1} f^0_s(g)=d_H(s)f^0_s(g)$$
by Lemma 6.39 as desired, completing the proof. 
(Compare with Lemma 5.2 of [GL].)

We note that by Lemma 6.43 the Mellin transform
$$\bL^{\std}_p (s)(m,g)=\int\limits_{M_{ab}(F)} |c(m,g)|^{s-\frac 12}\xi^0_P(m^{-1}g)\  dm$$
of $\xi^0_P$ is the basic function in the sense of doubling method as in equation (6.15) whose definition is justified by equation (6.17).

We now show that $\xi^0_P$ is preserved by Fourier transform $\sF_{\oP|P}$. Equation (5.8) of [BK2] can be stated as
$$\sF_{\oP|P} (\delta_{P,\gamma})={d_{\oP}\over d_P}\delta_{\oP,\gamma},\tag6.44$$
where $1/d_P :=d^{-1}_P$. We now apply (6.44) to $\xi^0_P=d^{-1}_{\oP}(\delta_{P,0})$ to get
$$\sF_{\oP|P}(\xi^0_P)=\xi^0_{\oP},$$
where $\xi^0_{\oP}=d^{-1}_P(\delta_{\oP,0})$. We have therefore proved:

\medskip\noindent
(6.45)
\ {\bf Proposition}
{\it The Fourier transform $\sF_{\oP|P}$ preserves our basic function $\xi^0_P$. More precisely,} 
$$\sF_{\oP|P}(\xi^0_P)=\xi^0_{\oP}.$$

\medskip\noindent
(6.46)\ {\bf Remark}.
We refer to the proof of Lemma 3.14 of [BK3] in Pages 12 and 13, where the calculations similar to our Lemma 6.43 and Proposition 6.45 are carried out.

\medskip\noindent
(6.47)\ {\bf Remark}.
By Part 3 of Theorem 5.10 of [BK2] and our discussion above, our basic function is in fact their function $c_{_{P,0}}$; also compare with formula (3.18) in Lemma 3.14, as well as Theorem 3.13 of [BK3]. But from [BK2] and [BK3], it was not clear if $c_{_{P,0}}$ is our basic function in the sense of giving the unramified $L$--function. 
Discussions with Wen--Wei Li after the first version of this manuscript was distributed, has led to an appendix [Li2] by him which also proves the equality of $c_{P,0}$ with our basic function $\bL^{\std}_P (s)$ up to a shift in $s$.
Finally, the last equation of [Li2] shows that as ``half--densities'' [BK2,Li], $c_P=c_{_{P,0}}
$ matches $\bL^{\std}_P(\frac 12)$, thus unifying [BK2,BNS,Li] with our results coming from the doubling method.

\parskip=0pt

\newpage
\noindent
Department of Mathematics

\noindent
Purdue University

\noindent
150 N.~University Street

\noindent
West Lafayette, IN\ \ 47907

\noindent
email:\ \ shahidi\@math.purdue.edu

\newpage

\centerline{\bf Appendix: A Comparison of Basic Functions}
\bigskip
\centerline{Wen-Wei Li}


\bigskip

\abstract{We show that in the doubling construction of Piatetski-Shapiro--Rallis, the basic functions defined by Shahidi and Braverman--Kazhdan are the same up to an explicit shift. We also discuss the {\it{raison d'\^etre}} of this shift.}
\endabstract
\endtopmatter

\medskip\noindent
{\bf 1.\ Introduction}

\bigskip\noindent

Shahidi [Sh] studied the doubling method of Piatetski-Shapiro--Rallis [PSR, GPSR] from the perspective of Braverman--Kazhdan [BK2]. In the unramified setting, he considered the {\it basic function} $\bL^{\std}_P(s)$ defined through inverse Satake transform of the standard $L$-factor for classical groups; see for example [BK1, {\S 5.7}] or [Li17] for an introduction to basic functions. A more detailed account of the doubling method will be given in {\S 2}. The main result (Theorem 5.2) of this appendix to [Sh] is that, up to an explicit shift in $s$, the function $c_P = c_{P,0}$ of [BK2, p.548] coincides with $\bL^{\std}_P(s)$ for symplectic groups. In {\S 6}, we will also relate that shift to the ubiquitous $\frac 12$ shift in the doubling method; see (6.2).

I am deeply grateful to F.\ Shahidi for encouragements and numerous discussions on this topic. I also thank B.\ C.\ Ng\^o for his explanations on the reference [BK].

\bigskip
\noindent
{\bf Notations}\  For a non-archimedean local field $F$, we write $|\cdot|$ for its normalized absolute value, $\frak{o}_F$ for its ring of integers with maximal ideal $\frak{p}_F$, and $q := |\frak{o}_F/\frak{p}_F|$ for its residual cardinality.

The modulus character $\delta_\Gamma$ of a locally compact group $\Gamma$ is characterized by $\text{d}\mu(gxg^{-1}) = \delta_\Gamma(g) \text{d}\mu(x)$ for any left Haar measure $\mu$.

For a scheme $S$ over a commutative ring $A$, we write $S(A)$ for its set of $A$-points. We use $H_{\der}$ (resp.\ $H_{\text{ab}}$) to denote the derived subgroup (resp.\ abelianization) of a reductive group scheme $H$; the Langlands dual group of $H$ is denoted by $\hat{H}$.

Unless otherwise specified, algebraic groups act on the right of varieties. In particular, $\text{GL}(V)$ acts on the right of a space $V$; this is consistent with [Li].

\bigskip\noindent
{\bf 2.\ Background}
\bigskip\noindent
We begin by reviewing the doubling construction for symplectic groups. Let $F$ be a field with $\text{char}(F) \neq 2$. Let $G = \text{Sp}(V) \subset \text{GL}(V)$ where $(V, \langle \cdot | \cdot \rangle)$ is a $2n$-dimensional symplectic $F$-vector space. Put $V^\square := V \oplus V$ equipped with the symplectic form $\langle \cdot | \cdot \rangle \oplus - \langle \cdot | \cdot \rangle$, so there is an embedding
$$ G \times G \hookrightarrow G^\square := \text{Sp}(V^\square). $$
Let $P \subset G^\square$ be the Siegel parabolic subgroup stabilizing the diagonal image of $V$, which is a Lagrangian in $V^\square$; let $M \simeq \text{GL}(V)$ be its Levi quotient. As in [Sh] and [Li, \S 7.1], we have the spaces
$$ X_P := P_{\der} \backslash G^\square \hookrightarrow X := \overline{X_P}^{\text {aff}} \quad \text{(affine closure)}. $$
Let $M_{\text{ab}} \times G^\square$ act on the right of $X_P$ by
$$ P_{\der}x \overset{(m,g)}\to\longmapsto P_{\der} m^{-1}xg, \quad (m,g) \in M_{\ab} \times G^\square. $$
Also let $M_{\ab} \times G \times G$ act on the right of $M_{\ab} \times G$ by
$$	(t, g) \overset{(a ,g_1, g_2)}\to\longmapsto (ta, g_2^{-1} g g_1).\tag2.1$$
Then $M_{\ab} \times G$ embeds into $X_P$ as the open $M_{\ab} \times G \times G$-orbit $X^+ \subset X$, which contains the coset $P_{\der}$; see [Li, \S 7.2].

As shown in 
[Li, \S 7.4], $X$ is a {\it{normal reductive monoid}} with unit group $M_{\ab} \times G$. Its smooth locus is precisely $X_P$.

It is convenient to identify $M_{\ab} \times G$ with $\bG_m \times G$ via
$$ M_{\ab} \simeq \text{GL}(V)_{\ab}\overset \det^{-1}\to{\underset \sim\to\longrightarrow} \bG_m. $$
Define the homomorphism
$$	c: M_{\ab} \times G \to M_{\ab} \overset\det^{-1}\to\longrightarrow \bG_m.\tag2.2$$
This is analogous to [Sh, (6.2)]. Furthermore, $X$ is a {\it{flat monoid}} in Vinberg's sense [Vin], and $c$ is the restriction to $M_{\ab} \times G$ of the {\it{abelianization map}} $X \to \bG_a$, still denoted by $c$. The point of using $\det^{-1}$ is that $c(x) \to 0$ when $x$ approaches the boundary $X \smallsetminus X^+$; cf.\ [Li, Lemma 7.2.5].

\medskip
\noindent
{\it{Remark} }2.1 \ 
	The orthodox way for looking at monoids is to consider $M_{\ab} \times G \hookrightarrow X$ as an $(M_{\ab} \times G)^2$-equivariant map, by composing with $M_{\ab}^2 \to M_{\ab}$, $(a,b) \mapsto ab^{-1}$. This is irrelevant since $M_{\ab}$ is a torus.

Hereafter we will work in the unramified setting, so that $F$ will be a non-archimedean local field with odd residual cardinality $q$. The hyperspecial subgroups $G(\frak{o}_F)$, $M_{\ab}(\frak{o}_F)$, etc.\ are also chosen. Unless otherwise stated, the Haar measures are normalized so that the hyperspecial subgroups have volume $1$.

The {\it{basic function}} in [Sh, (6.16)] is a function $\bL^{\std}_P(s)$ on $(M_{\ab} \times G)(F)$ depending on a complex variable $s$; the definition is the same as that of [BK1] and [Li17], which we will recall later. On the other hand, Braverman and Kazhdan defined in [BK2, p.548] a function $c_P = c_{P,0}$ over $X_P(F)$. The aim of this note is to elucidate their relations.

\medskip
\noindent
{\it{Remark} }2.2. \ 
	We confine ourselves to the symplectic case in order to use the results of [Li, \S 7] safely. Nonetheless, a generalization to other classical groups seems within reach.

\bigskip\noindent
{\bf 3.\ Global models in equal characteristics}

\bigskip\noindent
Note that $X_P \subset X$ and $c_P$ are also defined in the equal-characteristic case $F = \bF_q(\!(t)\!)$. In this set-up, $X$ and $X_P$ actually come from $\bF_q$-schemes of finite type. Let $\sL X$ be the formal arc space (over $\bF_q$) of $X$, so that $X(\frak{o}_F) = X(\bF_q [\![ t ]\!]) =  \sL X(\bF_q)$; cf. [BNS].

Fix a prime number $\ell$ with $\ell \nmid q$. Choose $\sqrt{q}$ inside $\overline{\bQ_\ell}$. Then $c_P$ takes value in $\overline{\bQ_\ell}$ by its definition in [BK2].

One advantage of the equal-characteristic set-up is the existence of global models of the singularities for $\sL X$, which we now recall in greater generality.

Let $C$ be a smooth proper geometrically connected curve over $\bF_q$, and suppose that $F$ is the local field attached to $v \in C(\bF_q)$.

Assume $X$ to be a normal affine $\bF_q$-variety, on which a connected reductive $\bF_q$-group $H$ acts on the right with an open dense orbit $X^0$ contained in the smooth locus. We have the formal arc space $\sL X$ as before.

\bigskip
\noindent
{\bf Definition 3.1.}
	The space $\sM_{X,H}$ over $\bF_q$ of non-based {\it quasi-maps} into $X$ is the open substack of $\text{Map}(C, [X/H])$, algebraic and locally of finite type, that maps each test $\bF_q$-scheme $S$ to the groupoid of data
	$$ E: \;\text{an $H$-torsor over}\; C \times S, \quad \phi: C \times S \to X \wedge^H E $$
	such that $\phi^{-1}(X^0 \wedge^G E)$ is open in $C \times S$ and surjects onto $S$; i.e.\ $\phi$ ``generically'' lands in $X^0 \wedge^H E$. See [BNS, \S 2]. In a similar manner, we define $\sL^\circ X \subset \sL X$ consisting of formal arcs that ``generically'' lands in $X^0$, see {\it op.\ cit.} We have
	$$ \sL^\circ X(\bF_q) = X(\frak{o}_F) \cap X^0(F). $$

By incorporating a trivialization of $E$ over the formal completion along $v \times S$ into the data $(E, \phi)$, we have the $\sL H$-torsor $p: \tilde{\sM}_{X,H} \to \sM_{X,H}$. By the discussions in [BNS, \S 2], there is a diagram
$$\sL^\circ X \overset h\to\longleftarrow \tilde{\sM}_{X,H} \overset p\to\longrightarrow \sM_{X,H}. $$
Roughly speaking, one defines $h$ by using the trivialization $\xi$ to assign a point of $\sL X$ from $(E, \phi)$.

\bigskip
\noindent
{\bf Definition 3.2.}
	We say that $x \in \sL^\circ X(\bF_q)$ and $m \in \sM_{X,H}(\bF_q)$ are {\it related} if
	
{$\bullet$ \ } $\phi$ lands in the smooth locus of $X$ off $v$;

{$\bullet$ \ } there exists $\tilde{m}$ such that $p(\tilde{m})=m$ and $h(\tilde{m})=x$.

Now we invoke the theory [BK] of IC-functions on $\sL^\circ X(\bF_q)$ (i.e.\ alternating sum of traces of Frobenius at stalks). The convention here is that IC-functions take value $1$ over the smooth stratum; i.e.\ the IC-sheaf is normalized to be $\overline{\bQ_\ell}$ on the main stratum, see [BK, Proposition 8.9]. The same normalization is applied to IC-sheaves of $\sM_{X,H}$.

We denote the IC-functions by $\text{IC}_{\sL X}$, $\sM_{X,H}$, etc.

\medskip
\noindent
{\it{Remark} }3.3. \ 
	Note that in [BK], the IC-function is defined on a subspace $\sL^\bullet X$ which is independent of $H$-action and contains the $\sL^\circ X$ of [BNS].

We will need the following examples.

\item{$\bullet$\ }  When $H = M_{\ab} \times G^\square$ and $X = \overline{X_P}^{\text{aff}}$ as in the doubling construction, $\sM_{X,H} = \overline{\text{Bun}_P}$ is Drinfeld's compactification. In this case $X^0 = X_P$. See [Sa, 3.3.2] or [Br, 2.4] for further details.

\item{$\bullet$ \ }  When $H = G' \times G'$ and $X$ is a normal reductive monoid with unit group $G'$, this is the global model considered in [BNS, \S 2]. In this case $X^0 = G'$.
	

\bigskip
\noindent
{\bf Theorem 3.4.} \ (Local-global compatibility of IC-functions).
{\it  when $m$ and $x$ are related as in Definition 3.2, we have}
	$$ \text{IC}_{\sM_{X,H}}(m) = \text{IC}_{\sL X}(x). $$                                                                                  

\noindent
{\it Proof.}
	See \S 0 or Proposition 9.2 of [BK]. The case of monoids is already in [BNS, (2.5)], and the general case makes use of the arguments thereof.

Note that we do not require the full formalism in [BK] of perverse sheaves, duality, etc.\ on $\sL X$: only the IC-functions and the local-global compatibility matter.

\bigskip\noindent
{\bf 4.\ Inverse Satake transform}
\bigskip\noindent
Details of the materials below can be found in [Li17]. Denote the normalized valuation of $F$ as $v$, so that $|\cdot| = q^{-v(\cdot)}$. Let $X$ be a normal reductive monoid over $F$ which is flat, with unit group $G'$, and let $c: X \to \bG_a$ be the abelianization map. Suppose that $\xi$ is a $G'(\frak{o}_F)$-bi-invariant function on $G'(F)$ such that $\xi_n$, the restriction of $\xi$ to $G'(F)^{v \circ c = n}$, is compactly supported for all $n \geq 0$. Then we can extend the usual Satake transform to $\xi$ by setting
$$
	\text{Sat}(\xi) = \sum_{n \geq 0} \text{Sat}(\xi_n).\tag4.1$$
The formal sum above lives in some completion of the range of the usual $\text Sat$.

Suppose that $X = X_\rho$ is associated to an irreducible representation $\rho$ of $\widehat{G'}$, as explained in [Li17, BNS, Sh]. The {\it basic function} $\bL^\rho(s)$ (with $s \in \bC$) is defined by
$$\gather
	\bL^\rho(s) = \sum_{n \geq 0} \bL^\rho(s)_n \quad \text{as above}, \\
	\forall n \geq 0, \quad {\text Sat} (\bL^\rho(s)_n) = \text{Tr} \left(\Sym^n \rho\right) q^{-ns}.
\endgather$$
Hence
$$ {\text Sat}(\bL^\rho(s)) = \sum_{n \geq 0} \text{Tr} \left( \Sym^n \rho \right) q^{-ns}. $$
The right-hand side can be evaluated at an irreducible unramified representation $\Pi$ when $\Re(s) >\!>_\Pi 0$, by evaluating at $\text{Tr}(\Sym^n \rho)$ at the Satake parameter of $\Pi$ in the adjoint quotient $\widehat{G'} /\!/ \widehat{G'}$. This yields another characterization
$$ {\text Sat}(\bL^\rho(s))(\Pi) = L(s,\Pi,\rho). $$

The definition of basic functions leads to
$$ \bL^\rho(s+t) = |c|^t \cdot \bL^\rho(s), \quad s,t \in \bC. $$
See [Li17], the discussions after Definition 3.2, where one writes $f_{\rho, s}$ instead of $\bL^\rho(s)$.

We note that the Satake transform is actually defined over $\bZ[\sqrt{q}]$ in the natural bases. The structure of inverse Satake transform and basic functions is determined entirely by the root datum of $G'$, $c$ together with $q$.

\bigskip\noindent
{\bf 5.\ The basic function of Braverman--Kazhdan}

\bigskip\noindent
The conventions are as in \S 2. Identify $M_{\ab}$ with $\bG_m$ via $\det^{-1}$ as before. Let $\rho := \id \boxtimes \std$ be the standard representation of $(\bG_m \times G)^\wedge$, with $\bG_m^\wedge$ acting by dilation. It is irreducible.
\bigskip\noindent
{\bf 5.1\ proposition}\ 
	{\it The representation $\rho$ corresponds to the $L$-{\text monoid} $X$ by the recipe in [BNS, \S 4].}

\medskip\noindent
{\it Proof.}\ 	This is [Li, Theorem 7.4.9].

Fix a prime number $\ell$ with $\ell \nmid q$, and fix $\iota: \overline{\bQ_\ell} \simeq \bC$ to reconcile with the harmonic analysis in [Sh]. The dual groups $\hat{G}$, etc.\ are taken over $\overline{\bQ_\ell}$. Take $\sqrt{q}$ inside $\overline{\bQ_\ell}$ via $\iota$. Then $c_P$ takes value in $\overline{\bQ_\ell}$.

It is clear from the formulas in [BK2] that $c_P$ depends only on

\item{$\bullet$ \ }   the residual cardinality $q$ of $F$,

\item{$\bullet$ \ }  some representation theory of the dual groups.

\medskip\noindent
Strictly speaking, the conventions in [BK2] are slightly different from ours: see the Remark 5.3.

The following result will be proved by passing to global models.
\bigskip\noindent
{\bf Theorem 5.2.}
{\it The following properties hold for $c_P$.
\item{1.\ } As a function on $X_P(F)$, it is $(M_{\ab} \times G^\square)(\frak{o}_F)$-invariant, and supported on $X(\frak{o}_F)$.
\item{2.\ } Restrict $c_P$ to the unit group $(M_{\ab} \times G)(F)$. Then the extended Satake transform (4.1) for $M_{\ab} \times G$ is applicable to $c_P$.
\item{3.\ } More precisely, we have
		$$ \bL^\rho(-n) = c_P. $$}

{\it{Proof.}}\ 
	The invariance under $(M_{\ab} \times G^\square)(\frak{o}_F)$ is contained in [BK2, p.547]. To see that ${\text Supp}(c_P) \subset X(\frak{o}_F)$, one invokes the Cartan decomposition for $X(\frak{o}_F)$ in [Sa, Theorem 2.3.8].
	
	The remaining assertions are of a combinatorial nature, depending solely on the residual cardinality $q$; see the discussions at the end of \S 4. Thus we can and do switch to the case $F = \bF_q(\!(t)\!)$. Let $C := \bP^1_{\bF_q}$ and choose $v \in C(\bF_q)$ to perform local-global arguments.

	In [Sa, pp.647--648] Sakellaridis affirms that
	$$ c_P = \Phi^0 $$
where $\Phi^0$ is the IC-function of $\overline{\text{Bun}_P}$ normalized as before, evaluated at $x \in X(\frak{o}_F) \cap X_P(F)$ by relating it to suitable points $m \in \overline{\text{Bun}_P}(\bF_q)$ as in Definition 3.2. We refer to {\it op.\ cit.} for details. The displayed equality is in turn a consequence of [BFGM, Theorem 7.3].

	Next, the local-global compatibility (Theorem 3.4) of IC-functions entails that
	$$ c_P = \Phi^0 = \text{IC}_{\sL X} $$
	as functions on $X(\frak{o}_F) \cap X_P(F) = \sL^\circ X(\bF_q)$.

	Recall that $\text{IC}_{\sL X}$ depends only on $X$ and not on the groups acting on it. We may restrict the $M_{\ab} \times G^\square$-action to $M_{\ab} \times G \times G$ (or inflate to $(M_{\ab} \times G)^2$ as in Remark 2.1), and regard $X$ as a normal reductive monoid with unit group $M_{\ab} \times G$. This operation shrinks $\sL^\circ X$ as $X^0$ is shrunk from $X_P$ to $M_{\ab} \times G$; we restrict $c_P$ and $\text{IC}_{\sL X}$ accordingly. Cf.\ Remark 3.3.
	
	In this setting of monoids, [BNS, Theorem 4.1], its errata [BNSe] together with Proposition 5.1 say that $\text{Sat}(c_P) = \text{Sat}(\text{IC}_{\sL X})$ is indeed well-defined by the recipe (4.1); furthermore, for all unramified irreducible representation $\chi \boxtimes \pi$ of $(M_{\ab} \times G)(F)$ we have 
	$$ \text{Sat}\left( c_P \right) \left( |c|^s \chi \boxtimes \pi \right) = L\left( s - \langle \eta_{M_{\ab} \times G}, \lambda \rangle, \chi \boxtimes \pi, \rho \right). $$
	Here $\Re(s) \gg_{\pi, \chi} 0$, and we choose a Borel subgroup $B \subset G$ to define
	
\item{$\bullet$ \ }   $\eta_{M_{\ab} \times G}$: the half-sum of $B$-positive roots of $M_{\ab} \times G$;

\item{$\bullet$ \ }   $\lambda$: the $\hat{B}$-highest weight of the representation $\rho$ of $(M_{\ab} \times G)^\wedge$.

\noindent
It remains to show $\langle \eta_{M_{\ab}\times G}, \lambda \rangle = n$. Let $B \twoheadrightarrow T$ be the Levi quotient. Choose the standard basis $\check{\epsilon}_1, \ldots, \check{\epsilon}_n$ for $X^*(T)$ such that the $B$-simple roots are
	 $$\check{\epsilon}_1 - \check{\epsilon}_2, \ldots, \check{\epsilon}_{n-1} - \check{\epsilon}_n, 2\check{\epsilon}_n. $$
	Let $\epsilon_1, \ldots, \epsilon_n$ be the dual basis for $X_*(T)$. It is clear that $M_{\ab}$ does not contribute to $\eta_{M_{\ab} \times G}$, thus
	$$ \eta_{M_{\ab} \times G} = \left(0, n \check{\epsilon}_1 + (n-1) \check{\epsilon}_2 + \cdots + \check{\epsilon}_n \right) \; \in \bZ \oplus X^*(T). $$
	By [Li, Proposition 7.4.8], $\lambda = (1, \epsilon_1) \in \bZ \oplus X_*(T)$. Our theorem follows.
\medskip
\noindent
Note that the proof involves two global models $\sM_{X,H}$ for the same local object $\sL X$: one with $H = M_{\ab} \times G^\square$ from [Sa], and the other with $H = M_{\ab} \times G \times G$ (or $(M_{\ab} \times G)^2$ by Remark 2.1) from [BNS].

\medskip\noindent
{\bf Remark}\  5.3.
	In [BK2] one considers $G^\square / P_{\text{der}}$ with left $G^\square$-action, and $M_{\ab}$ acts on it by $xP_{\text{der}} \overset m\to\mapsto xmP_{\text{der}}$. Here we follow the ``right'' conventions of \cite{Sake, Li}, thus the formula for $c_P$ must be modified. A formula for $c_P$ or $\Phi^0$ in our convention can be found in [Sa, (4-3)].

\bigskip\noindent
{\bf 6.\ Shifts}
\bigskip\noindent
Set $\rho := \id \boxtimes \std$ as before. The basic function $\bL_P^{\std}(s)$ in [Sh, (6.16)] is characterized by
$$ \text{Sat}(\bL_P^{\text std}(s))\left( |c|^s {\bold 1} \boxtimes \pi \right) = L(s, \bold 1 \boxtimes \pi, \rho) = L(s, \pi, \text{std}) $$
for all unramified irreducible representations $\pi$ of $G(F)$ and $\Re(s) \gg_\pi 0$. Since the unramified characters of $F^\times \simeq M_{\ab}(F)$ always take the form $|\cdot|^t$ for some $t \in \bC$, one can replace the trivial representation $\bold{1}$ by any unramified $\chi$. This implies
$$ \bL_P^{\text{std}}(s) = \bL^\rho(s)  = c_P |c|^{n+s}.$$

Shahidi [Sh] studies $\bL_P^{\text std}(s)$ from the viewpoint of doubling method. But what the doubling zeta integrals yield are $L(s + \frac{1}{2}, \ldots)$, which implies that $\bL^{\text std}_P(\frac{1}{2})$ is ``more basic'' than $\bL^{\text{std}}_P(0)$ from the doubling perspective. This is responsible for many occurrences of $s - \frac{1}{2}$ in [Sh]; cf.\ the Introduction in that reference. Hence we rewrite the equation above at $s=0$ as
$$
	\bL^{\text std}_P\left( \frac{1}{2} \right) = c_P \cdot |c|^{n + \frac{1}{2}}.\tag6.1$$

\noindent{\bf Question}. How to explain the shift in (6.1), admitting that $\bL^{\text std}_P(\frac{1}{2})$ and $c_P$ are both natural objects?

An analogous issue in Godement--Jacquet theory has been addressed at the end of the local part of [Li, \S 1.2], by using Schwartz--Bruhat {\it half-densities}, i.e.\ square-roots of measures. We set out to explicate the shift $n + \frac{1}{2}$ in (6.1) under the same paradigm, combining geometry and considerations from harmonic analysis.
\bigskip\noindent
{\bf Lemma 6.1.}\ 
	{\it The modulus character of $P(F)$ is}
	$$ \delta_P(m) = |\det(m)|^{-2n-1}, \quad m \in \text{GL}(V) \simeq M(F) \twoheadleftarrow P(F). $$

\noindent
{\it Proof.}
	If we assume that $\text{Sp}(V^\square) \subset \text{GL}(V^\square)$ acts \text{on the left } of $V^\square$, then $\delta_P(m) = |\det(m)|^{2n+1}$ is a well-known fact. One can pass to right actions by $gv = vg^{-1}$; this turns $\det$ into $\det^{-1}$, whence the assertion.

\medskip\noindent
{\bf Proposition 6.2.}\ 
{\it Regard $\delta_P$ as a function on $(M_{\ab} \times G)(F)$ via the projection to $M_{\ab}$. Then $\delta_P = |c|^{2n+1}$.}
\bigskip
\noindent
{\it Proof.}
	Combine Lemma 6.1 and (2.2).
\medskip
Note that $P_{\text{der}}(F)$ is (algebraically) unimodular. By choosing a nonzero $G^\square$-invariant algebraic volume form, we obtain a $G^\square(F)$-invariant measure $|\Xi|$ on $X_P(F)$. Likewise, denote a chosen Haar measure on $(M_{\ab} \times G)(F)$ by $|\Omega|$.

According to the $L^2$-philosophy in [Li] that is already manifest in [BK2], what is ``basic'' for $X_P$ is not the function $c_P$ but the half-density
$$ c_P |\Xi|^{\frac{1}{2}} $$
on $X_P(F)$. Likewise, the ``basic half-density'' for the monoid $X$ with unit group $X^+ \simeq M_{\ab} \times G$ is $\bL_P^{\std}\left( \frac{1}{2} \right) |\Omega|^{\frac{1}{2}}$ instead of $\bL_P^{\std}\left( \frac{1}{2} \right)$.

Denote the restrictions of $c_P$ and $|\Xi|$ to the open dense subspace $X^+(F)$ by the same symbols. To obtain an $(M_{\ab} \times G \times G)(F)$-invariant density on $X^+(F)$, by the discussions preceding [Li, Proposition 7.3.2], or a direct verification using the $M_{\ab}$-equivariance of  (2.1), we can take
$$ |\Omega| :=  \delta_P^{-1} |\Xi| $$
where $\delta_P$ is pulled-back to $(M_{\ab} \times G)(F)$. Hence
$$\align
	c_P |\Xi|^{\frac{1}{2}} & = c_P \delta_P^{\frac{1}{2}} \cdot |\Omega|^{\frac{1}{2}} \\
	& = c_P |c|^{n + \frac{1}{2}} \cdot |\Omega|^{\frac{1}{2}}
\endalign$$
by Proposition 6.2. 

We rescale $|\Xi|$ so that $|\Omega|$ gives mass $1$ to the hyperspecial subgroup of $(M_{\ab} \times G)(F)$. Then the formula above shows that (6.1) is equivalent to
$$	c_P |\Xi|^{\frac{1}{2}} = \bL^{\text{std}}_P\left( \frac{1}{2} \right) |\Omega|^{\frac{1}{2}},\tag6.2$$
an identity between two {\it basic half-densities} restricted to $X^+(F)$. This explains the $n + \frac{1}{2}$ shift in (6.1).

\parskip=0pt


\bigskip\bigskip
\noindent
E-mail address:  wwli\@math.ac.cn

\noindent
Academy of Mathematics and Systems Science, Chinese Academy of Sciences

\noindent
55, Zhongguancun donglu, 100190 Beijing, People's Republic of China.

\bigskip\noindent
University of Chinese Academy of Sciences 

\noindent
19A, Yuquan lu, 100049 Beijing, People's Republic of China.

\bye